\documentclass[aos,preprint,reqno]{imsart}

\setattribute{journal}{name}{}

\usepackage[utf8]{inputenc} 
\usepackage[T1]{fontenc}    
\usepackage[margin=1.5in]{geometry}
\usepackage{hyperref}       
\usepackage{url}            
\usepackage{booktabs}       
\usepackage{amsfonts}       
\usepackage{nicefrac}       
\usepackage{microtype}      

\usepackage{natbib}
\usepackage{amsmath, amssymb, amsthm, bbm, dsfont, mathrsfs}
\usepackage{amsthm}
\usepackage{subfiles, subcaption}
\usepackage{color, verbatim, graphicx}
\usepackage{cleveref}
\usepackage{mathtools}
\usepackage{commath}
\usepackage{stmaryrd}
\usepackage[symbol]{footmisc}

\newcommand{\FF}{\mathbb{F}}
\newcommand{\GG}{\mathbb{G}}

\newcommand{\XX}{\mathbb{X}}

\newcommand{\C}{\mathscr{C}}

\renewcommand{\P}{\mathscr{P}}
\newcommand{\KL}{\mathrm{KL}}

\DeclareMathOperator{\supp}{supp}

\newcommand{\nats}{\mathbb{N}}

\newcommand{\reals}{\mathbb{R}}
\newcommand{\spd}{\mathbb{S}_{++}}

\newtheorem{definition}{Definition}[section]
\newtheorem{theorem}[definition]{Theorem}
\newtheorem{proposition}[definition]{Proposition}
\newtheorem{lemma}[definition]{Lemma}
\newtheorem{corollary}[definition]{Corollary}

\newtheorem{assumption}[definition]{Assumption}

\begin{document}

\begin{frontmatter}
\title{Finite mixture models do not reliably learn \\ the number of components}
\runtitle{Finite mixtures do not reliably learn the number of components}

\begin{aug}
\author{\fnms{Diana} \snm{Cai}\thanksref{t1, paddr}},
\author{\fnms{Trevor} \snm{Campbell}\thanksref{t1, ubcaddr}},
\and
\author{\fnms{Tamara} \snm{Broderick}\thanksref{mitaddr}}

\thankstext{t1}{First authorship is shared jointly by D.\ Cai and T.\ Campbell.}

\affiliation{Princeton University\thanksmark{paddr}, University of British
    Columbia\thanksmark{ubcaddr}, \\ Massachusetts Institute of Technology\thanksmark{mitaddr}}

\address[paddr]{
    Department of Computer Science \\
    Princeton University \\
    Princeton, NJ, USA 08544\\
    Email: \href{mailto:dcai@cs.princeton.edu}{dcai@cs.princeton.edu}
}
\address[ubcaddr]{
    Department of Statistics \\
    University of British Columbia \\
    Vancouver BC V6T 1Z4 \\
    Email: \href{mailto:trevor@stat.ubc.ca}{trevor@stat.ubc.ca}
}
\address[mitaddr]{Computer Science and Artificial Intelligence Laboratory \\
   Massachusetts Institute of Technology \\
    Cambridge, MA, USA 02139
    Email: \href{tbroderick@csail.mit.edu}{tbroderick@csail.mit.edu}
}
\runauthor{D.\ Cai, T.\ Campbell, T.\ Broderick}
\end{aug}


\begin{abstract}
    Scientists and engineers are often interested in learning the number of subpopulations (or components) present in a data set. A common suggestion is to use a finite mixture model (FMM) with a prior on the number of components.  Past work has shown the resulting FMM component-count posterior is consistent; that is, the posterior concentrates on the true, generating number of components.  But consistency requires the assumption that the component likelihoods are perfectly specified, which is unrealistic in practice.  In this paper, we add rigor to data-analysis folk wisdom by proving that under even the slightest model misspecification, the FMM component-count posterior \emph{diverges}: the posterior probability of any particular finite number of components converges to 0 in the limit of infinite data.  Contrary to intuition, posterior-density consistency is not sufficient to establish this result.  We develop novel sufficient conditions that are more realistic and easily checkable than those common in the asymptotics literature. We illustrate practical consequences of our theory on simulated and real data.

\end{abstract}

\end{frontmatter}

\maketitle


\section{Introduction}
\label{sec-intro}

Mixture modeling is a mainstay of statistical machine learning. In applications where the number of mixture components
is unknown in advance, a principal inferential goal is to estimate and interpret this number.
For example, practitioners might wish to find the number of latent
genetic populations \citep{pritchard2000,lorenzen2006,huelsenbeck2007,tonkin2019fast},
gene tissue profiles \citep{yeung2001,medvedovic2002},
cell types \citep{chan2008,prabhakaran2016},
microscopy groups \citep{rubin2015bayesian,griffie2016bayesian},
haplotypes \citep{xing2006},
switching Markov regimes in US dollar exchange rate data \citep{otranto2002},
gamma-ray burst types \citep{mukherjee1998},
segmentation regions in an image (e.g., tissue types in an MRI scan \citep{banfield1993}),
observed humans in radar data \citep{teklehaymanot2018bayesian},
basketball shot selection groups \citep{hu2020bayesian},
or communities in a social network
\citep{geng2019probabilistic,legramanti2020extended}.

A natural question then is: can we reliably learn the number of latent groups in a data set?
To make this question concrete, we focus on a Bayesian approach.
Consider the case where the true, generating number of components is known.
A natural check on a Bayesian mixture analysis is to
establish that the Bayesian posterior on the number of components increasingly
concentrates near the truth as the number of data points becomes arbitrarily large.
In the remainder, we will focus on this check---though our work has practical implications beyond Bayesian
analysis.

A standard Bayesian analysis uses a component-count prior with support on all strictly-positive integers
\citep{miller2016}.
\citet{nobile1994} has shown that the component-count posterior of the resulting \emph{finite mixture model} (FMM) does
concentrate at the true number of components. But crucially,
this result depends on the assumption that
the component likelihoods are perfectly specified.
In every application we have listed above, the true generating component likelihoods
do not take a convenient parametric form that might be specified in advance.
Indeed, some form of misspecification, even if slight, is typical in practice.
So we must ask how the component-count posterior behaves
when the component likelihoods are misspecified.

\begin{figure*}[t]
    \centering
\begin{subfigure}[b]{0.46\linewidth}
    \centering
    \includegraphics[scale=0.42]{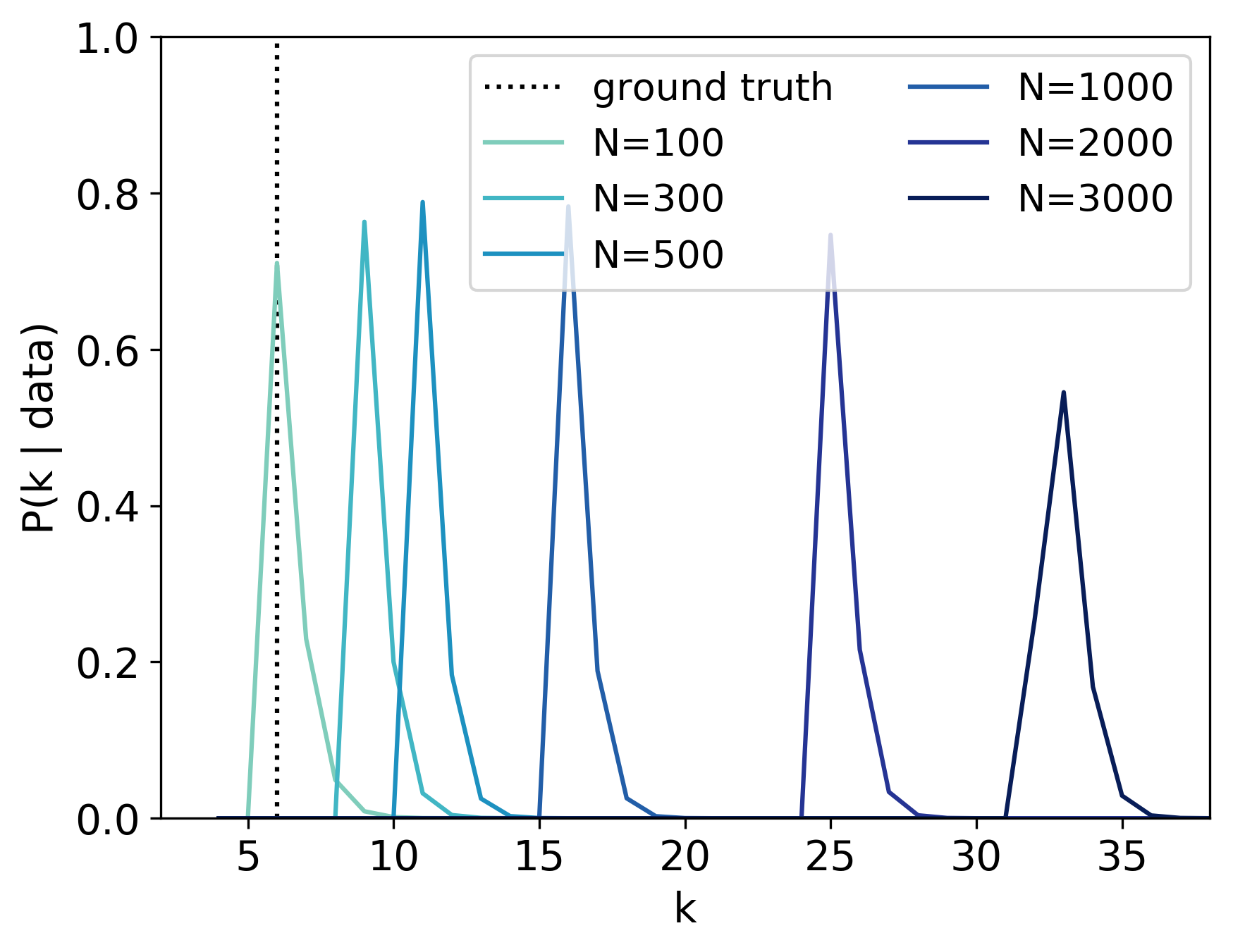}
    \caption{Mouse single-cell RNA-seq data}
    \label{fig-single-cell}
\end{subfigure}
\begin{subfigure}[b]{0.46\linewidth}
    \centering
    \includegraphics[scale=0.42]{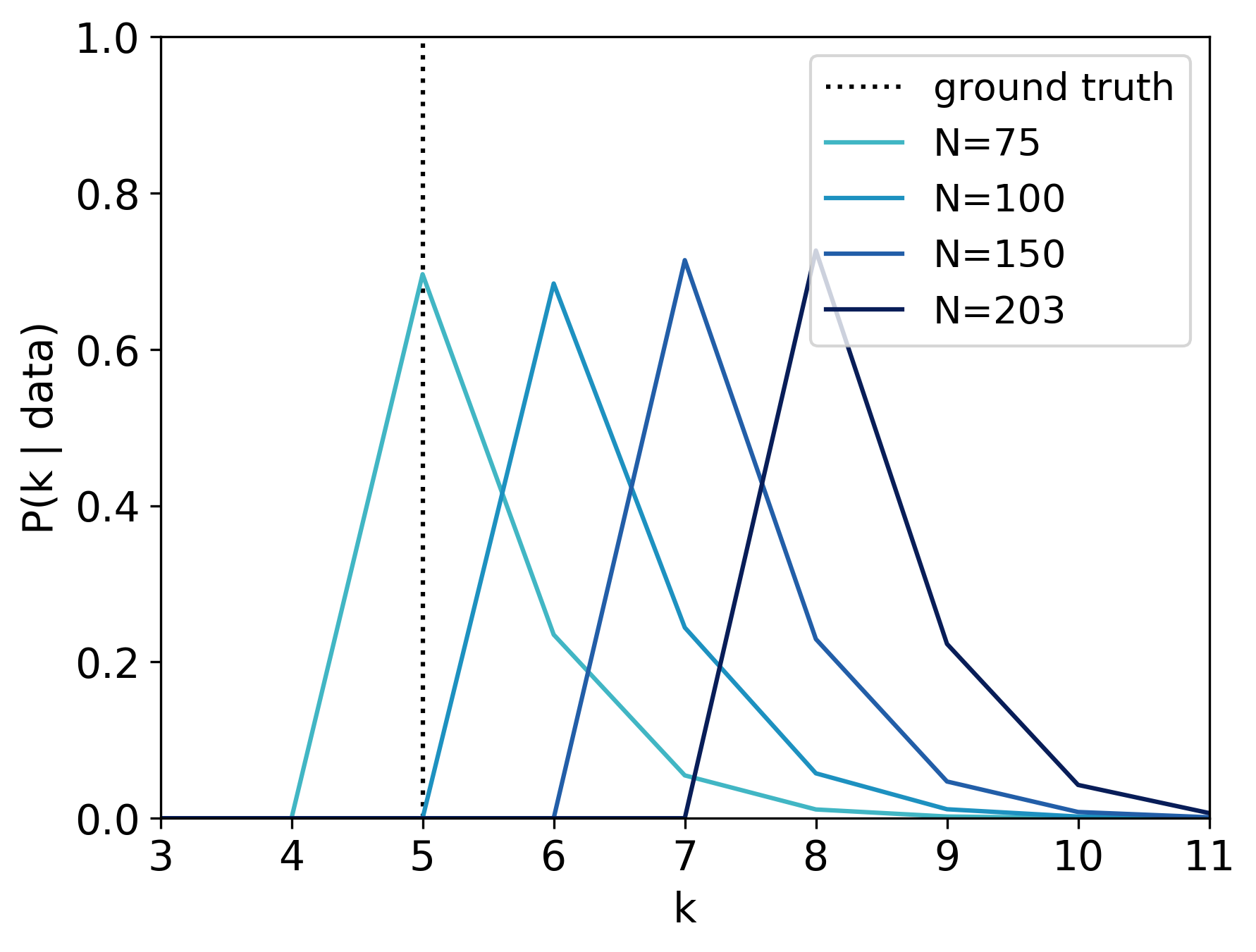}
    \caption{Lung tissue gene expression data}
    \label{fig-gene-expression}
\end{subfigure}
    \caption{
    Posterior probability of the number of components $k$ for Gaussian
    mixture models,
    fit to
    (a) mouse cortex single-cell RNA sequencing data
    and (b) lung tissue gene expression data.
    Details in \Cref{ssec-gene}.
    }
    \label{figure-data}
\end{figure*}

Data science folk wisdom suggests that when component likelihoods are
misspecified, mixture models will tend to overestimate the number
of clusters; see, e.g., Section 7.1 of \citet{fruhwirth2006finite}.
This overestimation is apparent in
\Cref{figure-data}, which shows the component-count posterior
of a Gaussian mixture model applied to two example gene expression data sets
\citep{de2008clustering,prabhakaran2016}.
In fact, \Cref{figure-data} demonstrates an effect far worse
than just overestimation:
the posterior distribution appears to concentrate for any (large enough) fixed
amount of data, but actually concentrates on increasing values as more data are observed.
Therefore, inference is unreliable;
the practitioner may draw quite different conclusions depending
on how large the data set is.

In the present paper, we add rigor to existing data science folk intuition
by proving that this behavior occurs in a wide class of FMMs under arbitrarily small
amounts of misspecification.
We examine FMMs with essentially any component shape---where we make only mild, realistic, and checkable
assumptions on the component likelihoods.
Notably, we include univariate and multivariate Gaussian component likelihoods
in our theory, but do not restrict only to these shapes.
We show that under our assumptions and when the component likelihoods are not perfectly specified,
the component-count posterior concentrates strictly away from the true number of components.
In fact, we go further to show that the FMM
posterior for the number of components \emph{diverges}: for \emph{any}
finite $k\in\nats$, the posterior probability that the number of components is
$k$ converges to 0 almost surely as the amount of data grows.

We start by introducing FMMs and
stating our main result in \Cref{sec-prelim}.
We discuss our assumptions in more detail in \Cref{sec-assumptions}
and prove our result in
\Cref{sec-proof}. In \Cref{sec-extensions} we extend our main theorem to priors that may vary as the data set grows.
We discuss related work below and in \Cref{sec-related}.
The paper concludes in \Cref{sec-experiments} with empirical evidence that the FMM component-count posterior
depends strongly on the amount of observed data.
Our results demonstrate that, in practice, past estimates of component number
may have strongly depended on the size of a particular data set.

\textbf{Filling a gap in the literature.}
While recent work has established various asymptotic properties of mixture models,
we observe that
our results here are not trivial extensions of existing research.
First note that, intuitively,
as the number of data points
grows, the posterior concentrates at the generating density \citep{schwartz1965,ghosh2003,ghosal2017fundamentals},
which can be well-approximated by an infinite mixture due in part to misspecification.
However, posterior consistency for the density alone
is not enough to guarantee consistency for the model parameters;
parameter consistency may not hold under, for instance, a discontinuous mapping from the
component parameter to
the component density.

Second, note that posterior divergence for the number of components could, in principle, be obtained
if parameter consistency for the mixture holds.
But existing results on parameter consistency, such as \citet{nguyen2013convergence},
focus on obtaining rates of contraction; thus these results rely on
stronger conditions that are typically verified for individual component families by imposing
additional constraints,
such as (second-order) strong identifiability of the mixture or a compact parameter space
\citep{chen1995optimal,nguyen2013convergence,heinrich2018strong}.
But neither of these constraints are satisfied by common families of interest such as Gaussians
with unknown mean and variance.
By contrast, our result uses the weakest notion of mixture identifiability
\citep{teicher1961identifiability}
along with a continuity condition on the component family,
and we relax the requirement of a compact parameter space.
To do so, we develop a novel theoretical condition that requires the component family to have \emph{degenerate
limits}. Together, these advances ensure
the applicability of our theory to practical likelihood families including, but not limited to, full Gaussians.
In fact,
the degenerate limits condition and its use in our analysis may be useful for extending other
results on posterior asymptotics that currently rely on compact parameter spaces.

Finally, \citet{miller13,miller2014} have shown that typical uses of Dirichlet process mixture models (DPMMs)
inconsistently estimate the true, generating number of components. But \citet{miller13,miller2014}
focus on the DPMM prior instead of the FMM and
on perfectly specified likelihoods. The DPMM is misspecified in a different sense than the one we focus on in the present paper:
namely, the DPMM uses infinitely many components though we assume finitely many generated the data. For this reason,
practitioners typically invoke the DPMM posterior on the number of \emph{clusters}
\citep{pella2006gibbs,huelsenbeck2007}, i.e., components represented in the observed data, rather than the component-count
posterior directly. Indeed, \citet{miller2016} recommend using the FMM we study here
to resolve the difficulties of the DPMM. Finally, observe that the work of
\citet{miller2016} demonstrates that nonparametrically estimating component shape with a DPMM would not
provide a simple resolution of the FMM divergence issue.


\section{Main result}
\label{sec-prelim}

We begin with a brief description of the finite mixture
model used in this work.
In this section, we provide just enough
detail to state \Cref{thm-inconsistency} and leave
the precise probabilistic details for \Cref{sec-assumptions}.
Let $g$ be a mixing measure
$g := \sum_{j=1}^k p_j \delta_{\theta_j}$
on a parameter space $\Theta$
with $p_j \in [0,1]$ and $\sum_{j=1}^k p_j=1$,
and let $\Psi = \{\psi_\theta: \theta \in \Theta\}$ be a family of component
distributions dominated by a $\sigma$-finite measure $\mu$.
We can express a finite mixture $f$
of the components as
\begin{align*}
      f &= \int_\Theta \psi_\theta dg(\theta) = \sum_{j=1}^k p_j \psi_{\theta_j}.
\end{align*}
Consider a Bayesian model with prior distribution $\Pi$ on the set of
all mixing measures $\GG$ on $\Theta$ with finitely many atoms, i.e.,
$g \sim \Pi$,
and likelihood corresponding to conditionally i.i.d.\ data from $f = \int \psi_\theta dg(\theta)$.
The model assumes the likelihood is
$f$, but the model is \emph{misspecified};
i.e.,
the observations $X_{1:N} := (X_1,\ldots,X_N)$
are  generated
conditionally i.i.d.\ from a finite mixture $f_0$ of
distributions \emph{not} in $\Psi$.

Our main result is that under this misspecification of the likelihood,
the posterior on the number of components
$\Pi(k \,|\, X_{1:N})$
diverges; i.e.,
for any finite $k\in\nats$, $\Pi(k\,|\, X_{1:N}) \to 0$ as $N\to\infty$.
We make only two requirements of the mixture model to guarantee this result:
(1) the true data-generating distribution $f_0$ must be arbitrarily well-approximated
by finite mixtures of $\Psi$, and (2) the family $\Psi$ must satisfy
mild regularity conditions that hold for popular mixture models
(e.g., the family $\Psi$ of Gaussians parametrized by mean and variance).
We provide precise definitions of the assumptions needed for \Cref{thm-inconsistency}
to hold in \Cref{sec-assumptions}, and a proof in \Cref{sec-proof}.

\begin{theorem}[Main result]
\label{thm-inconsistency}
    Suppose observations $X_{1:N}$ are
    generated i.i.d.~from a distribution $f_0$ that is not a finite mixture of $\Psi$. Assume that:
\begin{enumerate}
\vspace{-5pt}
\addtolength{\itemindent}{2cm}
    \item[\Cref{assumption-kl-support}:] $f_0$ is in the KL-support of the prior $\Pi$,
    \item[\Cref{assumption-regularity}:] $\Psi$ is continuous, is mixture-identifiable, and has degenerate limits.
\vspace{-5pt}
\end{enumerate}
    Then the posterior on the number of components
diverges; i.e., for all $k \in \nats$,
\begin{align}
    \label{eq-severe-inconsistency}
    \Pi(k \,|\, X_{1:N}) \stackrel{N\rightarrow\infty}{\longrightarrow} 0
 \quad f_0\text{-a.s.}
\end{align}
\end{theorem}

Note that the conditions of the theorem---although technical---are
satisfied by a wide class of models used in practice.
\Cref{assumption-kl-support} requires that the prior $\Pi$ places enough mass on
mixtures near the true generating distribution $f_0$.  \Cref{assumption-regularity} enforces regularity of the component family and is satisfied by many popular
models used in practice, such as the multivariate Gaussian family.
\begin{proposition}
    \label{prop-gaussians}
    Let
    $\Psi = \left\{ \mathcal{N}(\nu, \Sigma) : \nu \in \reals^d, \, \Sigma\in
    \spd^d \right\}$
    be the multivariate Gaussian family,
    where $\spd^d := \{\Sigma \in \reals^{d\times d}: \Sigma = \Sigma^\top, ~\Sigma \succ 0\}$
    is the set of $d \times d$ symmetric, positive definite matrices.
    Then $\Psi$ satisfies \Cref{assumption-regularity}.
\end{proposition}
Thus, provided that $f_0$ is in the KL-support of the prior,
under a misspecified Gaussian mixture model, our main result implies that the
posterior number of components
diverges.
While \Cref{prop-gaussians} is stated for Gaussian component distributions,
we generalize it to mixture-identifiable location-scale families $\Psi$ in \Cref{prop-locscale}.

Additionally, we note that the divergence of the posterior
given in \Cref{eq-severe-inconsistency}
is stronger than the behavior described in \citet{miller13}
for DPMMs: namely, \citet{miller13} show that
the posterior probability converges
to 0 at the \emph{true} number of components.
In contrast, here we show that the posterior
probability converges to 0 for any finite number of components.
We conjecture that posterior divergence also holds for DPMMs, but
the proof is outside of the scope of this paper.

\textbf{Extension: Priors that vary with $N$.} While the result of \Cref{thm-inconsistency} assumes that
the model uses a fixed prior $\Pi$, in practical modeling scenarios one may
specify a prior $\Pi_N$ that depends on the observed data $X_{1:N}$. For instance, these priors
can arise in empirical Bayes; see \Cref{sec-extensions,sec-experiments} for examples.
In \Cref{sec-extensions} we show that if $f_0$ satisfies a modified KL-support
condition with respect to the sequence of priors $\Pi_N$, the number of
components also diverges in this setting.

\textbf{Extension: Priors with an upper bound on the number of components.}
\Cref{thm-inconsistency} is designed for priors that place full support on any
positive integer number of components. One might instead use a prior that has
support on at most $\tilde{k}$ components, with $\tilde{k}$ finite. In this case,
the posterior number of components will not diverge to infinity but instead
typically concentrate on the upper bound, $\tilde{k}$. A precise statement of this behavior
appears in \Cref{thm-finite-prior} (in \Cref{appendix-finite-prior}) as an
analog of our \Cref{thm-inconsistency}.
\Cref{thm-finite-prior} shows that posterior inference does not improve over the
baseline estimate of the number of components provided by $\tilde{k}$. If
$\tilde{k}$ is already a good estimate of the number components, posterior
concentration at $\tilde{k}$ does not improve the estimate. In practice $\tilde{k}$
is often chosen as some large upper bound of convenience; then $\tilde{k}$ is not
a good estimate of the number of components, and concentration at $\tilde{k}$ is
undesirable.



\section{Precise setup and assumptions in \Cref{thm-inconsistency}}
\label{sec-assumptions}

This section makes the details of the modeling setup
and each of the conditions in
\Cref{thm-inconsistency} precise.

\subsection{Notation and setup}

Let $\XX$ and $\Theta$ be Polish spaces for the observations
and parameters, respectively, and endow both with their Borel $\sigma$-algebra.
For a topological space $(\cdot)$, let $\C(\cdot)$ be the bounded continuous functions from $(\cdot)$ into $\reals$,
and $\P(\cdot)$ be the set of probability measures on $(\cdot)$ endowed with the weak topology metrized by the
L\'evy-Prokhorov distance
$d$ (\Cref{def-levy-dist}).
We use $f_i \Rightarrow f$ and $f_i\iff f'_i$ to denote $\lim_{i\to\infty}d(f_i, f)=0$ and
$\lim_{i\to\infty} d(f_i, f'_i) = 0$, respectively, for $f_i, f'_i, f \in \P(\cdot)$.
We assume that the family of distributions $\Psi = \{\psi_\theta: \theta \in \Theta\}$
is absolutely continuous with respect to a $\sigma$-finite base measure $\mu$,
i.e., $\psi_\theta \ll \mu$ for all $\theta \in \Theta$,
and that for measurable $A\subseteq \XX$, $\psi_\theta(A)$ is a measurable function on $\Theta$.
Define the measurable
mapping $F : \P(\Theta) \to \P(\XX)$
from mixing measures to mixtures of $\Psi$, $F(g) = \int \psi_\theta dg(\theta)$.
Let $\GG$ be the set of atomic probability measures on $\Theta$ with finitely many atoms,
and let $\FF$ be the set of finite mixtures of $\Psi$.

In the Bayesian finite mixture model from \Cref{sec-prelim},
a mixing measure $g\sim \Pi$ is generated from a prior measure $\Pi$ on $\GG$,
and $f = F(g)$ is a likelihood distribution.

The posterior distribution on the mixing measure is,
for all measurable $A \subseteq \GG,$
\begin{align}
    \label{eq-posterior-distn}
    \Pi(A \,|\, X_{1:N}) =
    \frac{\int_A \prod_{n=1}^N \frac{df}{d\mu}(X_n)\, d\Pi(g)}
        {\int_{\GG} \prod_{n=1}^N \frac{df}{d\mu}(X_n) \, d\Pi(g)},
\end{align}
where $\frac{df}{d\mu}$ is the density of $f=F(g)$ with respect to $\mu$.
This posterior on the mixing measure $g\in\GG$ induces
a posterior on the number of components $k\in \nats$
by counting the number of atoms in $g$, and it also
induces a posterior on mixtures $f\in\FF$ via the
pushforward through the mapping $F$.
We overload the notation $\Pi(\cdot \,|\, X_{1:N})$ to refer to all of these posterior
distributions and $\Pi(\cdot)$ to refer to prior distributions; the meaning should be clear from context.

\subsection{Model assumptions}
The first assumption of \Cref{thm-inconsistency} is that while the true
data-generating distribution $f_0$ is not contained in the model class $f_0 \notin \FF$,
it lies on the boundary of the model class.  In particular, we
assume $f_0$ is in the \emph{KL-support} of the prior $\Pi$. Denote the
Kullback-Leibler (KL) divergence between probability measures $f_0$ and $f$ as
\[
\KL(f_0,f) := \left\{\begin{array}{ll}
\int \log\left(\frac{df_0}{df}\right) \, df_0 & f_0 \ll f\\
\infty & \text{otherwise}
\end{array}\right. .
\]
\begin{assumption}
    \label{assumption-kl-support}
   For all $\epsilon > 0$,
   the prior distribution $\Pi$ satisfies
    \begin{align*}
        \Pi(f\in\FF : \mathrm{KL}(f_0, f) < \epsilon) > 0.
    \end{align*}
\end{assumption}
We use \Cref{assumption-kl-support} in the proof of \Cref{thm-inconsistency}
primarily to ensure that the Bayesian posterior is consistent for $f_0$.
Note that \Cref{assumption-kl-support} is fairly weak in practice. Intuitively,
it just requires that the family $\Psi$ is rich enough so that mixtures of $\Psi$
can approximate $f_0$ arbitrarily well, and that the prior $\Pi$ places
sufficient mass on those mixtures close to $f_0$.
For Bayesian mixture modeling,
\citet[Theorem~3]{ghosal1999},
\citet[Theorem~3.2]{tokdar2006posterior},
\citet[Theorem~2.3]{wu2008kullback},
and
\citet[Theorem~1]{petralia2012repulsive}
provide conditions needed to satisfy \Cref{assumption-kl-support}.

The second assumption of \Cref{thm-inconsistency} is that
the family of component distributions $\Psi$ is well-behaved.
This assumption has three stipulations. First,
the mapping $\theta \mapsto \psi_\theta$ must be continuous; this condition essentially asserts that similar
parameter values $\theta$ must result in similar component distributions $\psi_\theta$.
\begin{definition}
The family $\Psi$ is \emph{continuous} if the map $\theta \mapsto \psi_\theta$ is continuous.
\end{definition}

Second, the family $\Psi$ must be  \emph{mixture-identifiable}, which
guarantees that each mixture $f\in \FF$ is
associated with a unique mixing measure $G\in\GG$.
\begin{definition}[\citet{teicher1961identifiability,teicher1963identifiability}]
\label{defn-mixture-identifiable}
The family $\Psi$ is \emph{mixture-identifiable}
    if the mapping $F(g) = \int \psi_\theta dg(\theta)$ restricted to finite mixtures $F : \GG \to \FF$
is a bijection.
\end{definition}
In practice, one should always use
an identifiable mixture model for clustering; without identifiability,
the task of learning the number of components is ill posed.
And many models satisfy mixture-identifiability, such as
finite mixtures of
the multivariate Gaussian family \citep{yakowitz1968identifiability},
the Cauchy family \citep{yakowitz1968identifiability}, the gamma family
\citep{teicher1963identifiability},
the generalized logistic family,
the generalized Gumbel family,
the Weibull family, and von Mises family \citep[Theorem~3.3]{ho2016strong}.
A number of authors \citep[e.g.][]{chen1995optimal,ishwaran2001bayesian,nguyen2013convergence,ho2016strong,guha2019,heinrich2018strong}
appeal to stronger notions of identifiability for mixtures than
\Cref{defn-mixture-identifiable}.
But, to show posterior divergence in the present work,
we do not require conditions stronger than \Cref{defn-mixture-identifiable}.

The third stipulation---that the family $\Psi$
has \emph{degenerate limits}---guarantees that a ``poorly behaved'' sequence of parameters $(\theta_i)_{i \in \nats}$
creates a likewise ``poorly behaved''
 sequence of distributions $(\psi_{\theta_i})_{i \in \nats}$.
This condition allows us to rule out such sequences in the
proof of \Cref{thm-inconsistency}, and is the essential regularity condition to guarantee that a sequence of
finite mixtures of at most $k$ components cannot approximate $f_0$ arbitrarily closely.
\begin{definition}\label{defn-muwide}
    A sequence of distributions $(\psi_i)_{i=1}^\infty$
    is \emph{$\mu$-wide} if for any closed set $C$ such that $\mu(C) = 0$
    and any sequence of distributions $(\phi_i)_{i=1}^\infty$ such that $\psi_i\iff\phi_i$,
    \[
    \limsup_{i\to\infty}\phi_i(C) = 0.
    \]
\end{definition}
\begin{definition}\label{defn-degenerate-limits}
    The family $\Psi$ has \emph{degenerate limits} if
    for any tight, $\mu$-wide sequence $(\psi_{\theta_i})_{i\in\nats}$,
    we have that $(\theta_i)_{i\in\nats}$ is relatively compact.
\end{definition}
The contrapositive of \Cref{defn-degenerate-limits}
provides an intuitive explanation of the condition:
as $i \rightarrow \infty$,
for any sequence of parameters $\theta_i$ that eventually leaves every compact set $K\subseteq \Theta$,
either the $\psi_{\theta_i}$ become ``arbitrarily flat'' (not tight) or ``arbitrarily peaky'' (not $\mu$-wide).
For example, consider the family $\Psi$ of Gaussians on $\reals$ with Lebesgue measure $\mu$. If the variance of $\psi_{\theta_i}$ shrinks as $i$ grows,
the sequence of distributions converges weakly to a sequence of point masses (not dominated by the Lebesgue measure). If either the variance
or the mean diverges, the distributions flatten out and the sequence is not tight. We use the fact that these are the only two possibilities
when a sequence of parameters is poorly behaved (not relatively compact) in the proof of \Cref{thm-inconsistency}.

These three stipulations together yield \Cref{assumption-regularity}.
\begin{assumption}
\label{assumption-regularity}
The mixture component family $\Psi$ is continuous, is mixture-identifiable, and has degenerate limits.
\end{assumption}



\section{Proof of \Cref{thm-inconsistency}}
\label{sec-proof}

The proof  has two essential steps.
The first is to show that the Bayesian posterior is weakly consistent for the mixture
$f_0$; i.e., for any weak neighborhood $U$ of $f_0$ the sequence of posterior distributions satisfies
\begin{align}
\Pi(U \,|\, X_{1:N}) \stackrel{N\rightarrow\infty}{\longrightarrow} 1, \quad
    \text{$f_0$-a.s.}\label{eq:claimweakconst}
\end{align}
By
Schwartz's theorem (\Cref{thm-schwartz}),
weak consistency for $f_0$ is guaranteed directly by \Cref{assumption-kl-support}
and
the fact that $\Psi$ is dominated by a $\sigma$-finite measure $\mu$.
The second step is to show that for any finite $k\in\nats$, there exists a weak neighborhood $U$ of $f_0$
containing no mixtures of the family $\Psi$ with at most $k$ components.
Together, these steps show that the posterior probability
of the set of all $k$-component mixtures converges to 0 $f_0$-a.s.~as the amount of observed data grows.

We provide a proof of the second step. To begin, note that \Cref{assumption-kl-support} has two additional implications
about $f_0$ beyond \Cref{eq:claimweakconst}. First, $f_0$ must be absolutely continuous with respect to the dominating
measure $\mu$; if it were not, then there exists a measurable set $A$ such that $f_0(A) > 0$ and $\mu(A) = 0$. Since
$\mu$ dominates $\Psi$, any $f\in\FF$ satisfies $f(A) = 0$. Therefore $\KL(f_0,f) = \infty$, and the prior support
condition cannot hold. Second, it implies that $f_0$ can be arbitrarily well-approximated by finite mixtures under the weak
metric, i.e., there exists a sequence of finite mixtures $f_i \in \FF$, $i\in\nats$ such that
$f_i \Rightarrow f_0$ as $i\to\infty$. This holds because $\sqrt{\frac{1}{2}\KL(f_0, f)} \geq \text{TV}(f_0, f) \geq d(f_0, f)$.

Now suppose the contrary of the claim for the second step, i.e., that there exists
a sequence $(f_i)_{i=1}^\infty$ of mixtures of at most $k$ components from $\Psi$ such that $f_i\Rightarrow f_0$.
By mixture-identifiability, we have a sequence of mixing measures $g_i$ with at most $k$ atoms
such that $F(g_i) = f_i$.
Suppose first that the atoms of the sequence $(g_i)_{i\in\nats}$ either
stay in a compact set or have weights converging to 0.
More precisely, suppose there exists a compact set $K \subseteq \Theta$ such that
\begin{align}
g_i\left(\Theta \setminus K\right) \to 0.\label{eq:gito0}
\end{align}
Decompose each $g_i = g_{i,K} + g_{i,\Theta\setminus K}$ such that
$g_{i,K}$ is supported on $K$ and $g_{i,\Theta\setminus K}$ is supported on $\Theta\setminus K$.
Define the sequence of probability measures $\hat g_{i,K} = \frac{g_{i,K}}{g_{i,K}(\Theta)}$ for sufficiently large $i$
such that the denominator is nonzero.
Then \Cref{eq:gito0} implies
\[
 F\left(\hat g_{i,K}\right)\Rightarrow f_0.
\]
Since $\Psi$ is continuous and mixture-identifiable, the restriction of $F$ to the domain $\GG$
is continuous and invertible; and
since $K$ is compact, the elements of $(\hat g_{i,K})_{i\in\nats}$
are contained in a compact set $\GG_K\subseteq \GG$ by Prokhorov's theorem (\Cref{thm-prokhorov}).
Therefore
$F(\GG_K) = \FF_K$ is also compact,
and the map $F$ restricted to the domain $\GG_K$
is uniformly continuous with a uniformly continuous inverse by \citet[Theorems~4.14, 4.17, 4.19]{rudin1976principles}.
Next since $F(\hat g_{i,K})\Rightarrow f_0$, the sequence $F(\hat g_{i,K})$ is Cauchy in $\FF_K$; and since $F^{-1}$ is uniformly
continuous on $\FF_K$, the sequence $\hat g_{i,K}$ must also be Cauchy in $\GG_K$.
Since $\GG_K$ is compact, $\hat g_{i,K}$ converges in $\GG_K$.
\Cref{lemma-k-comps} below guarantees that the convergent limit $g_K$ is also a mixing measure with at most $k$ atoms;
continuity of $F$ implies that $F(g_K) = f_0$, which is a contradiction, since by assumption $f_0$ is not representable as a finite mixture of $\Psi$.
\begin{lemma}\label{lemma-k-comps}
Suppose $\phi, (\phi_i)_{i\in\nats}$ are Borel probability measures on a Polish space
such that $\phi_i\Rightarrow \phi$ and $\sup_i |\supp \phi_i| \leq k\in\nats$. Then  $|\supp \phi|\leq k$.
\end{lemma}
\begin{proof}
Suppose $|\supp \phi| > k$. Then we can find $k+1$ distinct
points $x_1, \dots, x_{k+1} \in\supp\phi$. Pick any metric $\rho$ on the Polish space,
and denote the minimum pairwise distance between the points $2\epsilon$.
Then for each point $j=1, \dots, k+1$ define the bounded, continuous function
    $h_j(x) = 0\vee \left(1-\epsilon^{-1}\rho(x, x_j)\right)$.
Since $x_j \in \supp\phi$, we have that $\int h_jd\phi > 0$.
Weak convergence $\phi_i\Rightarrow\phi$ therefore implies $\min_{j=1,\dots,k+1}\liminf_{i\to\infty} \int h_jd\phi_i > 0$.
But the $h_j$ are nonzero on disjoint sets, and each $\phi_i$ only has $k$ atoms; the pigeonhole principle yields a contradiction.
\end{proof}

Now we consider the remaining case: for all compact sets $K\subseteq \Theta$, $g_i(\Theta \setminus K) \not\to 0$.
Therefore there exists a sequence of parameters $(\theta_i)_{i=1}^\infty$
that is not relatively compact such that $\limsup_{i\to\infty} g_i(\{\theta_i\}) > 0$.
By \Cref{assumption-regularity}, the sequence $(\psi_{\theta_i})_{i\in\nats}$ is either not tight
or not $\mu$-wide. If $(\psi_{\theta_i})_{i\in\nats}$ is not tight then  $f_i = F(g_i)$ is not tight,
and by Prokhorov's theorem $f_i$ cannot converge to a probability measure, which contradicts $f_i\Rightarrow f_0$.
If $(\psi_{\theta_i})_{i\in\nats}$ is not $\mu$-wide then $f_i = F(g_i)$ is not $\mu$-wide.
Denote $(\phi_i)_{i\in\nats}$ to be the singular sequence  associated with $(f_i)_{i\in\nats}$
and $C$ to be the closed set such that $\limsup_{i\to\infty} \phi_i(C) > 0$, $\mu(C) = 0$, and $\phi_i\iff f_i$ per \Cref{defn-muwide}.
Since $f_0 \ll \mu$, $f_0(C) = 0$.
But $f_i \Rightarrow f_0$ implies that $\phi_i\Rightarrow f_0$, so
 $\limsup_{i\to\infty} \phi_i(C) = f_0(C) = 0$ by the Portmanteau theorem
(\Cref{thm-portmanteau}).
 This is a contradiction.


\section{Extension to priors that vary with $N$}
\label{sec-extensions}

Our main result (i.e., \Cref{thm-inconsistency}) applies to the setting of a
fixed prior $\Pi$.
However, it is often natural to specify a prior distribution
that changes with $N$
(e.g., \citealp{roeder1997practical}; \citealp{richardson1997}; and \citealp[Section~7.2.1]{miller2016}).
\Cref{cor-priors} below demonstrates that a result nearly identical
to \Cref{thm-inconsistency} holds for priors that are allowed to
vary with $N$, provided that $f_0$ is in the KL-support of the
\emph{sequence} of priors $\Pi_N$. The only difference is that
our result in this case is slightly weaker: we show that the
 posterior number of components diverges in probability rather
than almost surely.

\begin{assumption}
\label{assumption-kl-sequence}
For all $\epsilon > 0$,
the sequence of prior distributions $\Pi_N$ satisfies
\begin{align*}
    \liminf_{N \rightarrow \infty} \Pi_N(f: \KL(f_0, f) < \epsilon) > 0.
\end{align*}
\end{assumption}

\begin{corollary}
\label{cor-priors}
Suppose in the setting of \Cref{thm-inconsistency}
we replace \Cref{assumption-kl-support} with \Cref{assumption-kl-sequence}.
Then the posterior on the number of components
diverges in $f_0$-probability: i.e., for all $k \in \nats$,
\begin{align*}
    \Pi(k \,|\, X_{1:N}) \stackrel{N\rightarrow\infty}{\longrightarrow} 0
    \quad \text{in~}f_0\text{-probability}.
\end{align*}
\end{corollary}

\begin{proof}
    Since  for any $\epsilon > 0$,
    $\liminf_{N \rightarrow \infty} \Pi_N(f: \KL(f_0, f) < \epsilon) > 0$,
    \citet[Theorem~6.17, Lemma~6.26, and Example~6.20]{ghosal2017fundamentals}
    imply that
    the posterior is weakly consistent at $f_0$ in probability:
    i.e., for any weak neighborhood $U$ of $f_0$,
    \begin{align*}
        \Pi(U \,|\, X_{1:N}) \stackrel{N \rightarrow\infty}{\longrightarrow} 1
        \quad \text{in~}f_0\text{-probability}.
    \end{align*}
\Cref{assumption-kl-sequence} also implies that
    for sufficiently large $N$,
    $f_0$ is a weak limit of finite mixtures in $\FF$. The remainder of the proof
    is identical to that of \Cref{thm-inconsistency}.
\end{proof}


\section{Related work}
\label{sec-related}

In this work, we consider FMMs with a prior on
the number of components.
In the broader Bayesian mixture modeling literature,
posterior consistency for the mixture density \citep{ghosal1999,lijoi2004,kruijer2010adaptive}
and the mixing measure \citep{nguyen2013convergence,ho2016strong,guha2019}
is well established.
But posterior consistency for the number of components is not as thoroughly characterized.
There are several results establishing consistency for the number of components in well-specified FMMs.
\citet[Proposition~3.5]{nobile1994} and \citet[Theorem~3.1a]{guha2019}
demonstrate that FMMs exhibit posterior consistency
for the number of components when the model is well specified and
$\Psi$ is mixture-identifiable.
The present work characterizes the behavior of the FMM posterior
on the number of components under component misspecification.
Under misspecification of the component family or the support of the true mixing measure,
\citet[Theorem~4.1, Theorem~4.3]{guha2019}
establish posterior rates of contraction for the mixing measure for
Gaussian and Laplace location mixtures with compact parameter spaces.
Our results, which rely on posterior density consistency results, assume
weaker conditions on the prior and hold for more general classes of component families,
such as multivariate Gaussians parameterized by a mean and covariance.

A related approach for handling a finite but unknown number of components
is to specify a prior with a finite upper bound on the number of components
\citep[e.g.][]{ishwaran2001bayesian,chambaz2008bounds,rousseau2011,malsiner2016model,zhang2018identification,fruhwirth2019here}.
In the setting of overfitted FMMs with well-specified component densities,
\citet[Theorem~1]{rousseau2011} show that under a stronger identifiability
condition than mixture-identifiability and additional regularity assumptions on the model,
the posterior will concentrate properly by emptying the extra components.
\citet[Theorem~1]{ishwaran2001bayesian}
consider the setting of estimating the number of components with the assumption of
a known upper bound on the true number of components and well-specified components, and
show that the posterior does not asymptotically underestimate the number of components
when assuming a stronger identifiability condition than mixture-identifiability
and a KL-support condition on the prior.
Under a weaker notion of (second-order) strong identifiability \citep{chen1995optimal} and
a well-specified model,
\citet[Theorem~4]{chambaz2008bounds} provide upper bounds on the underestimation
and overestimation error of the number of components; furthermore, they
show that their conditions are satisfied by univariate Gaussians with bounded mean and variance
\citep[Corollary~1]{chambaz2008bounds}.
Notably, all of these methods with finite-support priors assume well-specified component families. By contrast, we show in Theorem B.1 that even for these finite-support priors, misspecified component families yield unreliable estimates of the number of components.

\citet{fruhwirth2006finite} provides a wide-ranging review of methodology for finite mixture modeling. In (e.g.) Section 7.1, \citet{fruhwirth2006finite} observes that, in practice, the learned number of mixture components will generally be higher than the true generating number of components when the likelihood is misspecified---but does not prove a result about the number of components under misspecification.
Similarly, \citet[Section~7.1.5]{miller2016} discuss the issue of estimating
the number of components in FMMs under model misspecification and state
that the posterior number of components is expected to diverge to infinity as
the number of samples increases, but no proof of this asymptotic behavior is
provided.

Finally, a growing body of work is focused on developing more robust FMMs and related mixture models. In order to address the issue of component misspecification, a number of authors propose using finite mixture models with nonparametric component densities, e.g. Gaussian-mixture components \citep{bartolucci2005clustering,di2007mixture,malsiner2017identifying} or overfitted-mixture components \citep{aragam2018identifiability}. However, for these finite mixture models that have mixtures as components, the posterior number of components and its asymptotic behavior have yet to be characterized.



\section{Experiments}
\label{sec-experiments}

\begin{figure*}[t!]
\centering
\begin{subfigure}[b]{0.46\linewidth}
    \centering
    \includegraphics[scale=0.40]{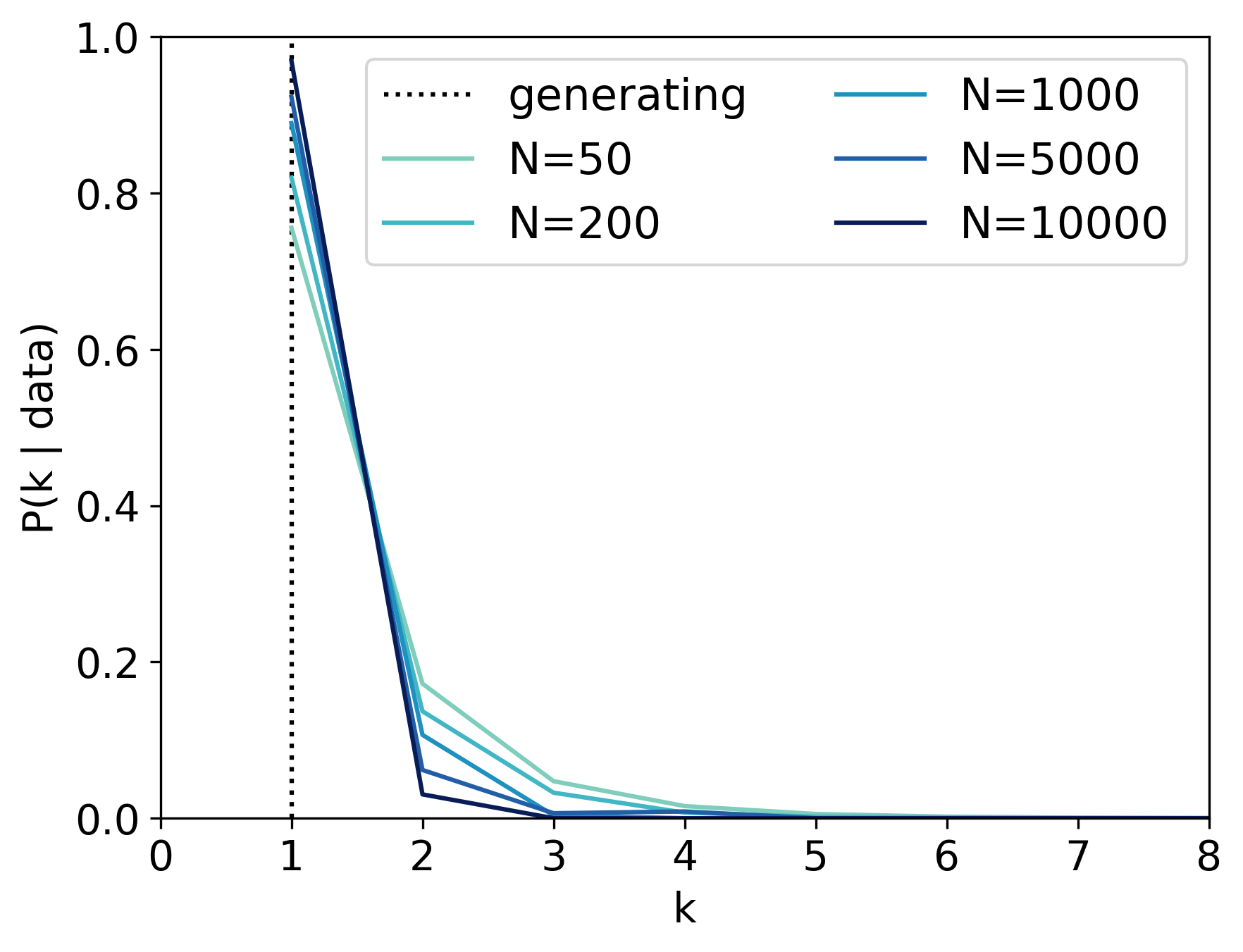}
    \caption{Gaussian data, 1 component}
    \label{fig-1comp-gaussian}
\end{subfigure}
\begin{subfigure}[b]{0.46\linewidth}
    \centering
    \includegraphics[scale=0.40]{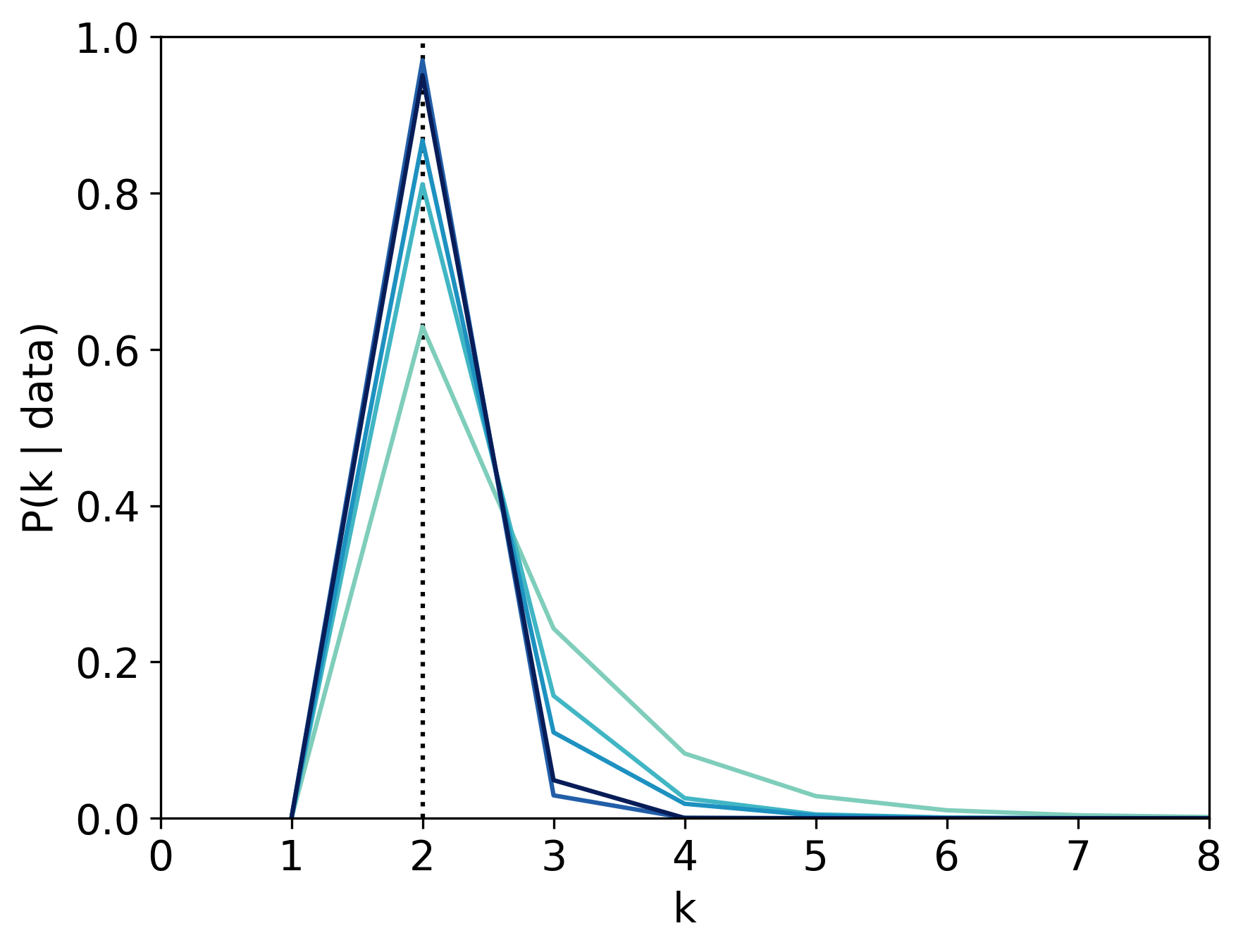}
    \caption{Gaussian data, 2 components}
    \label{fig-2comp-gaussian}
\end{subfigure}
\begin{subfigure}[b]{0.46\linewidth}
    \centering
    \includegraphics[scale=0.40]{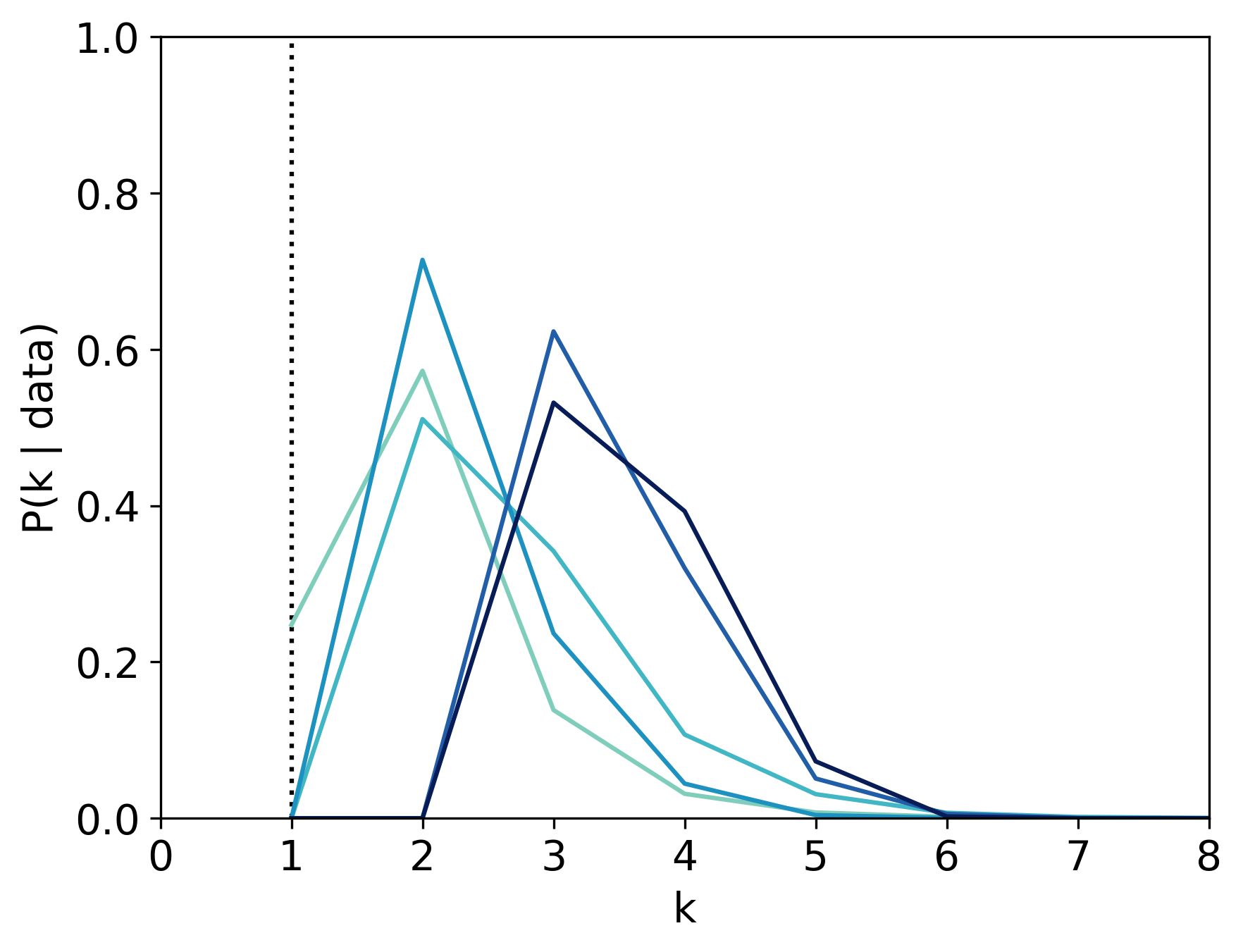}
    \caption{Laplace data, 1 component}
    \label{fig-1comp-laplace}
\end{subfigure}
\begin{subfigure}[b]{0.46\linewidth}
    \centering
    \includegraphics[scale=0.40]{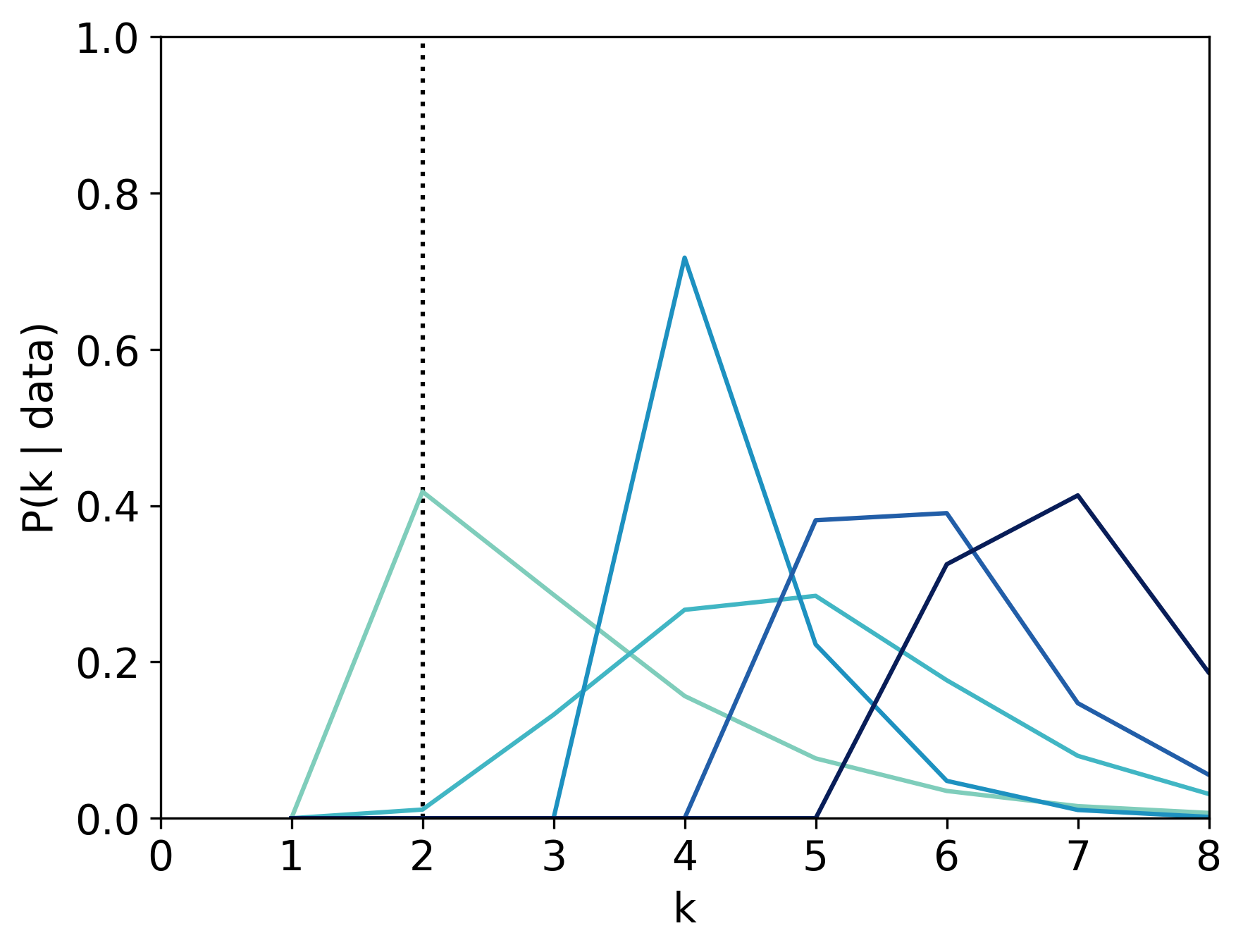}
    \caption{Laplace data, 2 components}
    \label{fig-2comp-laplace}
\end{subfigure}
\begin{subfigure}[b]{0.46\linewidth}
    \centering
    \includegraphics[scale=0.40]{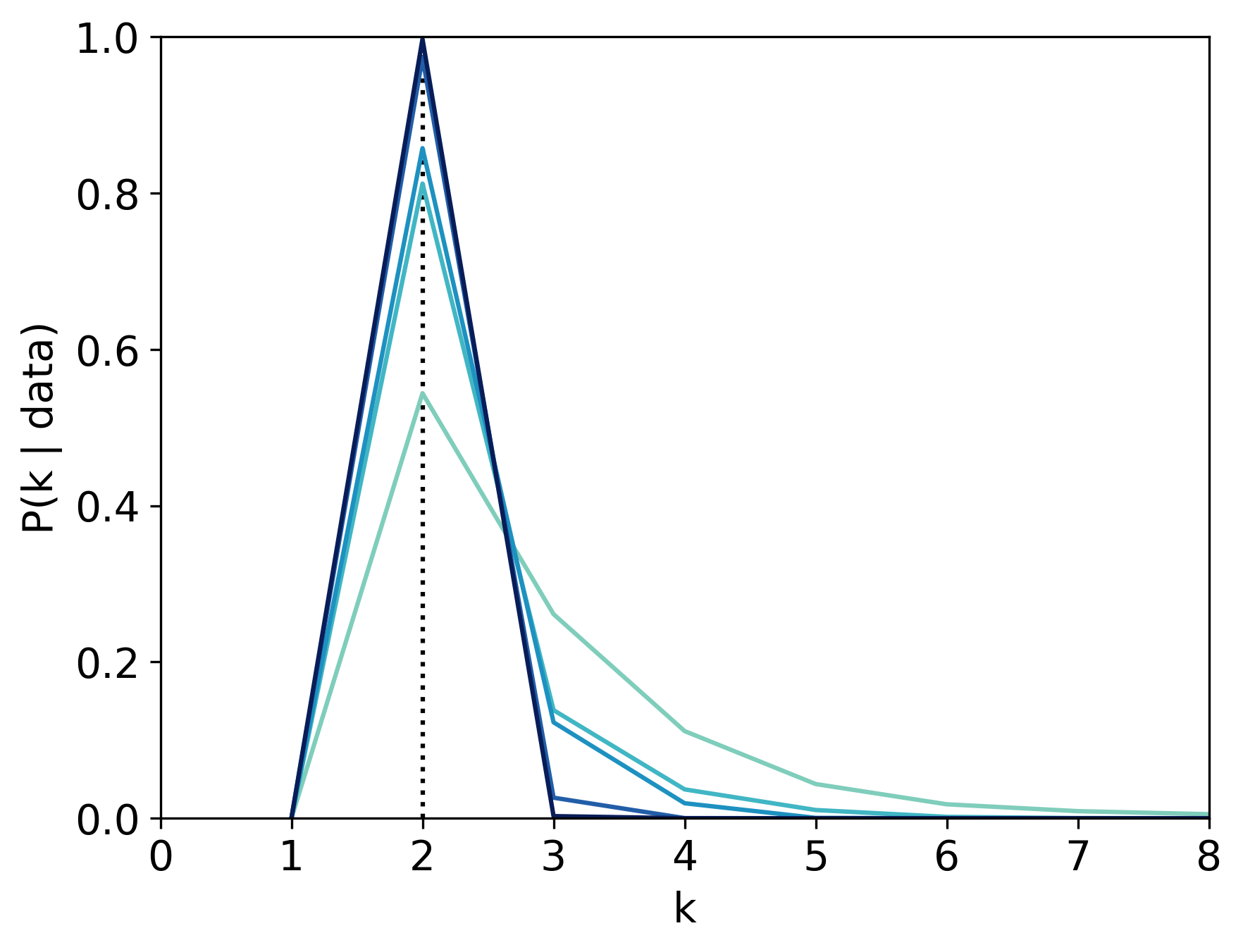}
    \caption{Gaussian data, varying prior}
    \label{fig-gaussian-fixed}
\end{subfigure}
\begin{subfigure}[b]{0.46\linewidth}
    \centering
    \includegraphics[scale=0.40]{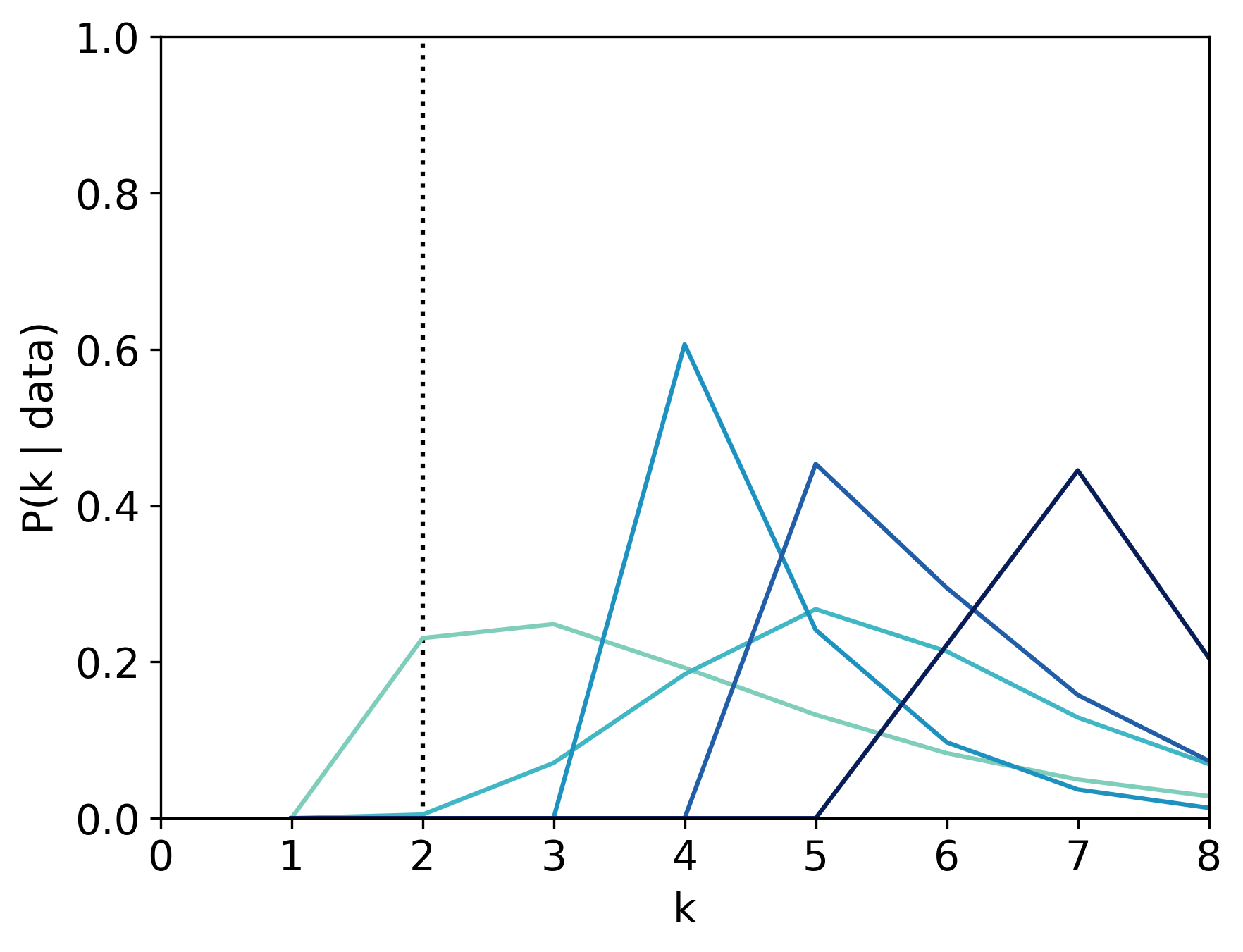}
    \caption{Laplace data, varying prior}
    \label{fig-laplace-varying}
\end{subfigure}
    \caption{
    \emph{Upper and middle rows}: Posterior probability of the number of components $k$ for Gaussian
    mixture models with a fixed prior
    fit to (a,b)
    univariate data generated from a Gaussian
    mixture model
    and (c,d) a Laplace mixture model,
    \emph{Lower row}: Posterior probability of the number of components of Gaussian mixtures with a varying prior fit to
    (e)
    2-component univariate data from a Gaussian mixture model and
    (f)
    2-component univariate data from a Laplace mixture model.
    }
    \label{fig-sims}
\end{figure*}

In this section, we demonstrate one of the primary practical implications
of our theory: the inferred number of components can change
drastically depending on the amount of observed data in misspecified
finite mixture models.
For all experiments below, we use a finite mixture model with
a multivariate Gaussian component family having diagonal covariance matrices
and a conjugate prior on each dimension.
In particular, consider number of components $k$, mixture weights $p\in\reals^k$,
Gaussian component precisions $\tau \in \reals_+^{k\times D}$ and means $\theta\in\reals^{k\times D}$,
labels $Z\in\{1, \dots, k\}^N$, and data $X\in\reals^{N \times D}$.

Then the probabilistic
generative model is
\begin{align*}
k &\sim \mathrm{Geom}(r) & p &\sim \mathrm{Dirichlet}_k(\gamma, \dots, \gamma)\\
\tau_{jd} &\overset{\text{i.i.d.}}{\sim} \mathrm{Gam}(\alpha, \beta)& \theta_{jd}
    &\overset{\text{i.i.d.}}{\sim} \mathcal{N}(m, \kappa_{jd}^{-1})\\
Z_{n}&\overset{\text{i.i.d.}}{\sim} \mathrm{Categorical}(p) & X_{nd}&\overset{\text{ind}}{\sim} \mathcal{N}(\theta_{z_nd}, \tau^{-1}_{z_nd}),
\end{align*}
where $j$ ranges from $1,\dots, k$, $d$ ranges from $1,\dots,D$, and $n$ ranges from $1,\dots,N$.
For posterior inference, we use a Julia implementation of
 split-merge collapsed Gibbs sampling \citep{neal2000markov,jain2004split}
from \citet{miller2016}.\footnote{Code available at
https://github.com/jwmi/BayesianMixtures.jl.}
The model and inference algorithm are described
in more detail in \citet[Sec.~7.2.2, Algorithm~1]{miller2016}.
Note that we use this model primarily
to illustrate the problem of posterior divergence under model misspecification;
it should not be interpreted as a carefully-specified model
for the data examples that we study.
Also note that while the empirical examples below involve Gaussian FMMs,
our theory applies to a more general class of component distributions.

\subsection{Synthetic data}
\label{ssec-synthetic}
\paragraph{Gaussian and Laplace mixtures}
Our first experiments on synthetic data are inspired by Figure~3 of \citet{miller2015},
which investigates the posterior of a mixture of perturbed Gaussians.
Here we study the effects of varying
data set sizes under both well-specified and misspecified models.
We generated data sets of increasing size $N \in \{50,200,1000,5000,10000\}$ from
1- and 2-component univariate Gaussian and Laplace mixture models,
where  the 1-component distributions have
mean 0 and scale 1,
and
the 2-component distributions have
means $(-5, 5)$, scales $(1.5, 1)$,
and mixing weights $(0.4, 0.6)$.
We generated the sequence of data sets such that each was a subset of the next, larger data set in the sequence.
Following \citet[Section~7.2.1]{miller2016},
we set the hyperparameters of the Bayesian finite mixture model as follows:
$m = \frac{1}{2}(\max_{n \in [\tilde N]}X_n + \min_{n \in [\tilde N]}X_n)$ where $\tilde N = 10{,}000$,
$\kappa = (\max_{n \in [\tilde N]}X_n  - \min_{n \in [\tilde N]}X_n)^{-2}$,
$\alpha=2$,
$r = 0.1$, $\gamma = 1$,
and $\beta \sim \mathrm{Gam}(0.2, 10/\kappa)$.
We refer to \citet[Section~7.2.1]{miller2016}
for additional details on the choice of model
hyperparameters and the sampling of $\beta$.
We ran a total of 100,000 Markov chain Monte Carlo iterations per data set; we
discarded the first 10,000 iterations as burn-in.

The results of the simulations are shown in
\Cref{fig-sims}.
For the data generated from the 1-component models,
the posterior on the number of components concentrates around 1
in the case of Gaussian-generated data as the sample size
increases (\Cref{fig-1comp-gaussian}),
whereas the posterior on the number of components
diverges for the Laplace data (\Cref{fig-1comp-laplace}).
We observe similar behavior in the
2-component case, where the posterior concentrates around the correct value in
the Gaussian case (\Cref{fig-2comp-gaussian})
but not the Laplace case (\Cref{fig-2comp-laplace}).

\paragraph{Priors that vary with $N$}
Next, we considered the same finite Gaussian mixture model described above but with a prior that varies with the data.
Specifically, for the prior on the means, we set the hyperparameters to
$m_N = \frac{1}{2}(\max_{n \in [N]}X_n + \min_{n \in [N]} X_n)$ and
$\kappa_N = (\max_{n \in [N]}X_n  - \min_{n \in [N]} X_n)^{-2}$,
which is the setting considered by \citet[Section~7.2.1]{miller2016};
the other hyperparameters were otherwise set to the same values above.
We used the 2-component Gaussian and Laplace data sets constructed above for the
fixed prior case.
The bottom row of \Cref{fig-sims} shows the results of the posterior
number of components under this prior for the well-specified and
misspecified cases; again we observe that the posterior diverges
under model misspecification.

\begin{figure*}[t]
\centering
\begin{subfigure}[b]{0.46\linewidth}
    \centering
    \includegraphics[scale=0.42]{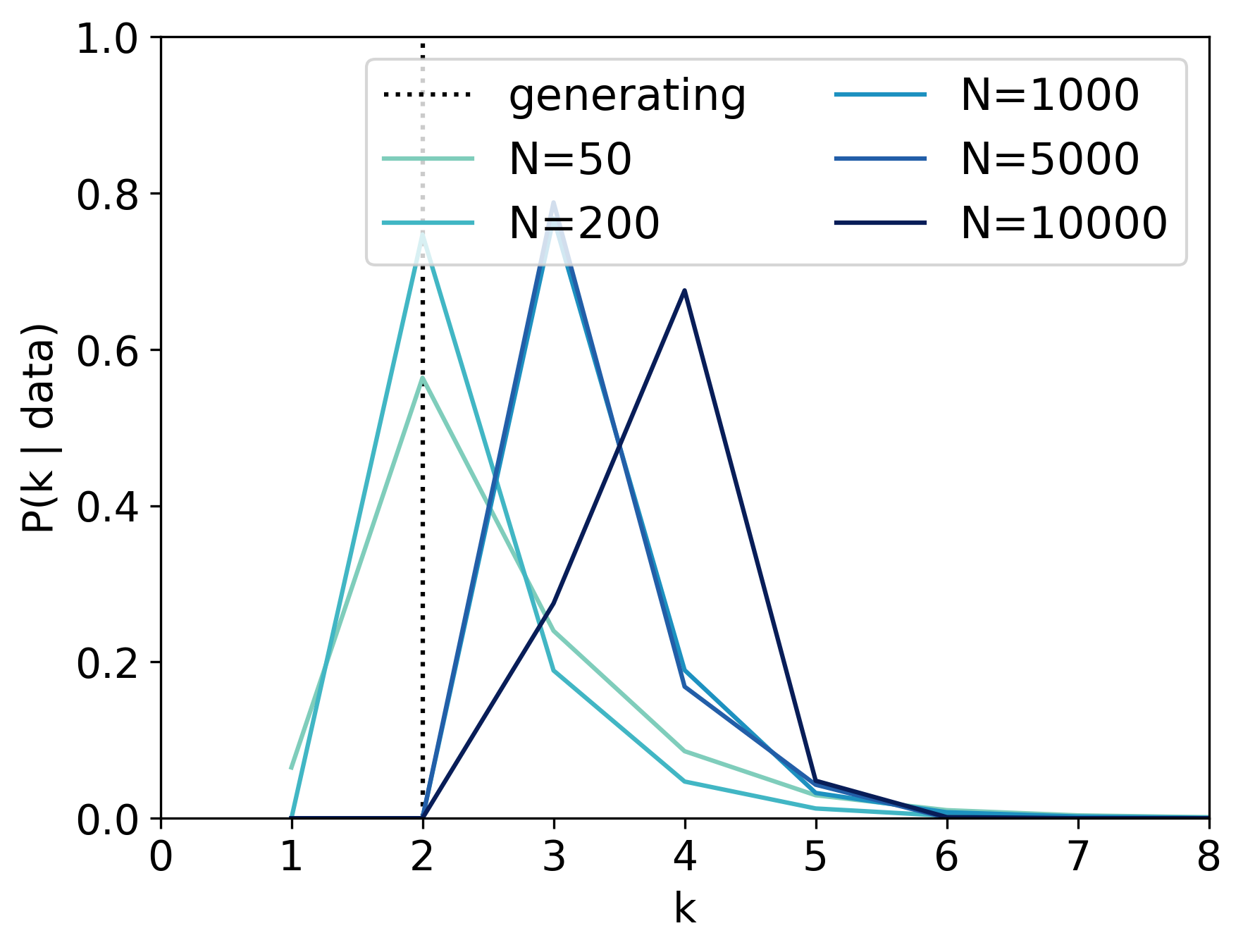}
    \caption{Data generated with $\epsilon=0.01$}
    \label{fig-contamin1}
\end{subfigure}
\begin{subfigure}[b]{0.46\linewidth}
    \centering
    \includegraphics[scale=0.42]{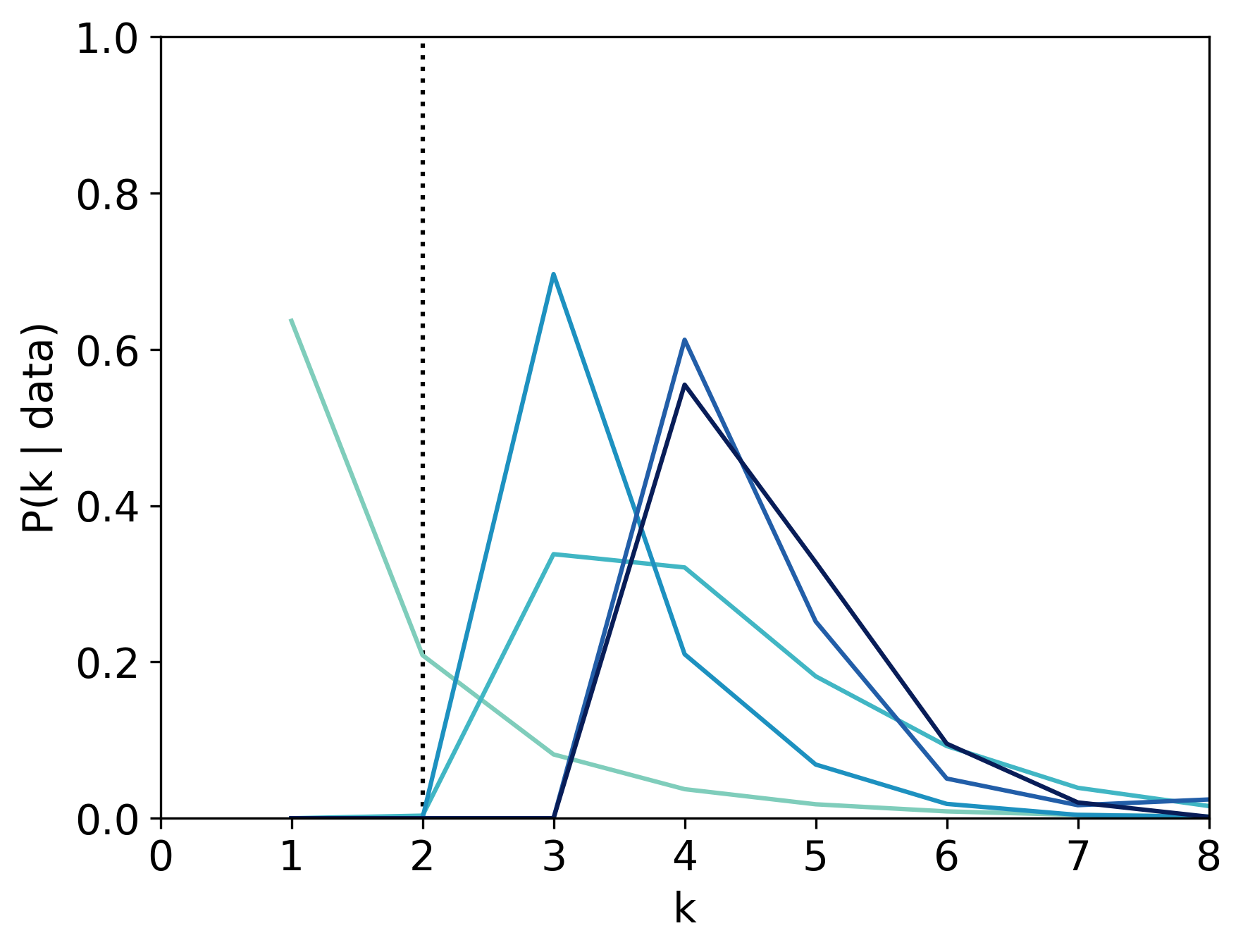}
    \caption{Data generated with $\epsilon=0.1$}
    \label{fig-contamin2}
\end{subfigure}
\caption{
    Posterior probability of the number of components $k$ for Gaussian
    mixture models with a fixed prior fit to data
    generated from an $\epsilon$-contaminated 2-component Gaussian mixture model,
    where
    $\epsilon$ is the proportion of data generated from a Laplace distribution.
}
    \label{figure-contamination}
\end{figure*}

\paragraph{$\epsilon$-contamination}
Finally, in order to study the posterior number of components
under a very slight amount of misspecification,
we applied the fixed-prior Gaussian mixture model above
to data generated with $\epsilon$-contamination.
That is, we generated the data according to the $\epsilon$-contaminated distribution
    $f_0 = (1-\epsilon) f  + \epsilon q,$
where $f$ is a 2-component Gaussian mixture distribution
with means $(5,10)$, variances $(1, 1.5)$, and mixing weights $(0.4,0.6)$,
and
$q$ is a Laplace distribution with location 0 and scale 1.
We generated two data sets: one with $\epsilon = 0.01$ and one with $\epsilon = 0.1$.
In \Cref{figure-contamination},
we observe that even under very small amounts of misspecification,
the posterior number of components diverges.

\subsection{Gene expression data}
\label{ssec-gene}

Computational biologists are interested in classifying cell types by applying
clustering techniques to gene expression
data \citep{yeung2001,medvedovic2002,mclachlan2002mixture,medvedovic2004bayesian,rasmussen2008modeling,de2008clustering,mcnicholas2010model}.
In our next set of experiments, we apply the Gaussian finite mixture model
to two gene expression data sets:
(1) single-cell RNA sequencing data from mouse cortex and hippocampus cells \citep{zeisel2015cell}
with the same feature selection as \citet{prabhakaran2016}
($N=3008$, $D=558$, 11,000 Gibbs sampling steps with 1,000 of those as burn-in)
and (2) mRNA expression data from human lung tissue \citep{bhattacharjee2001classification}
($N=203, D=1543$, and 10,000 Gibbs sampling steps with 1,000 of those burn-in).
Our experiments here represent a simplified version of previous mixture model analyses
for these and other related data sets \citep{de2008clustering,prabhakaran2016,armstrong2001mll,miller2016}.

As these gene expression data sets contain counts, we first transformed
the data to real numerical values. In particular, we used a base-2 log transform
followed by standardization---such that each dimension of the data had
zero mean and unit variance---per standard practices (e.g., \citet{miller2016}).
Then to examine the effect of increasing data set size on inferential results,
we randomly sampled subsets of increasing size without replacement;
each smaller subset was contained in the next larger data set.
For both data sets, we used hyperparameters $\alpha=1$, $\beta=1$, $m = 0$,
$\kappa_{jd} = \tau_{jd}$,
$r = 0.1$, and $\gamma = 1$.

For the single-cell RNAseq data set, the posterior on the number of
components is shown in \Cref{fig-single-cell}.
Here the ground truth number of clusters
is captured when the data set size is $N=100$.
But as predicted by our theory,
as we increase the number of data points,
the posterior number of components diverges.

The posterior on the number of components for the lung gene expression data
is shown in \Cref{fig-gene-expression}.
Again we find that on the smallest data subsets,
the posterior appears to capture the ground truth number of clusters,
but that as we examine more and more data, the posterior diverges.

The diagonal covariance Gaussian components are a particularly simple form of cluster shape.
But no matter how complex the component model, one could wonder whether an
even-more complex model might solve the
issue that the number of components diverge. In the typical real-world situation that the
component model cannot be specified in absolute perfection,
our theory confirms that the divergence problem will remain.
Thus, these examples suggest the need for more robust analyses.



\section{Discussion}
\label{sec-discussion}

We have shown that the posterior distribution
for the number of components in finite mixtures
diverges when the mixture component family is misspecified.
Since misspecification is almost unavoidable in real applications,
it follows that finite mixture models are typically unreliable
for estimating the number of components.
In practice, our conclusion implies that inferences on the number of
components can change drastically depending on the size of the data set, calling
into question the usefulness of these counts in application.

Since our analysis is inherently asymptotic, it is
possible that the Bayesian component-count posterior may still
provide useful inferences for a finite sample---for instance if care is taken
to account for the dependence of inferential conclusions on data set size.
Indeed, a number of authors have recently proposed robust Bayesian inference methods
to mitigate likelihood misspecification
\citep{woo2006robust,woo2007robust,rodriguez2011nonparametric,grunwald2014, miller2015,
bissiri2016general,
wang2017, holmes2017assigning, jewson2018principles, huggins2019,
knoblauch2019generalized,rigon2020generalized};
it remains to better understand connections between our results and these
methods.


\subsection*{Acknowledgments}
We thank Jeff Miller for helpful conversations and comments on a
draft.
D.\ Cai was supported in part by a Google Ph.D.\ Fellowship in Machine Learning.
T.\ Campbell was supported by a National Sciences and Engineering Research Council of Canada (NSERC) Discovery Grant and an NSERC Discovery Launch Supplement.
T.\ Broderick was supported in part by ONR grant N00014-17-1-2072, an
MIT Lincoln Laboratory Advanced Concepts Committee Award, a Google
Faculty Research Award, the CSAIL--MSR Trustworthy AI Initiative, and an ARO YIP award.

\appendix

\section{Background on posterior contraction}

\paragraph{Weak topology and weak convergence}

In this section, we review some definitions and results used in our work.
Our treatment follows
\citet[Appendix~A]{ghosal2017fundamentals}, and we refer to this chapter for additional details
on the topology of weak convergence.

Let $\XX$ be a Polish space metrized by $\rho$.
Below we define the L\'evy-Prokhorov metric, which induces the
weak topology on $\mathscr{P}(\mathbb{X})$.
\begin{definition}
\label{def-levy-dist}
Let $f, g \in \mathscr{P}(\mathbb{X})$.
The L\'evy-Prokhorov metric is defined as
    \begin{align*}
    d(f, g) =
    \inf\{
        \epsilon > 0: f(A) < g(A^\epsilon) + \epsilon,
        g(A) < f(A^\epsilon) + \epsilon
    \},
    \end{align*}
    where $A^\epsilon := \{y: \rho(x,y) < \epsilon \text{~for~some~} x \in A\}$.
\end{definition}

The Portmanteau theorem characterizes equivalent
notions of weak convergence, and
below we include the relevant portions of the
Portmanteau theorem used in our proof.
For a full statement of the theorem,
see \citet[Theorem~A.2]{ghosal2017fundamentals}.
\begin{theorem}[Portmanteau (partial statement)]
    \label{thm-portmanteau}
    The following statements are equivalent for any
    $f_i, f \in \mathscr{P}(\XX)$:
    \begin{enumerate}
        \vspace{-5pt}
        \item $f_i \Rightarrow f$;
        \item for all bounded, uniformly continuous
            $h: \mathbb{X} \rightarrow \reals$,
            \[\int h df_i \longrightarrow \int h df;\]
        \item
            for every closed subset $C$,
            $\limsup_i f_i(C) \leq f(C)$.
    \end{enumerate}
\end{theorem}

Prokhorov's theorem \citep[Theorem A.4]{ghosal2017fundamentals}
characterizes (weakly) compact subsets
of $\P(\XX)$ in terms of a \emph{tight} subset of measures.
A subset $\Gamma \subseteq \mathscr{P}(\mathbb{X})$ is
\emph{tight} if for any $\epsilon > 0$, there exists a compact subset
$K_\epsilon \subseteq \XX$ such that
for every $\psi \in \Gamma$,
$\psi(K_\epsilon) \geq 1 - \epsilon$.

\begin{theorem}[Prokhorov]
    \label{thm-prokhorov}
    If $\mathbb{X}$ is a Polish space,
    then $\Gamma \subseteq \mathscr{P}(\mathbb{X})$ is
    relatively compact if and only if $\Gamma$ is tight.
\end{theorem}

\paragraph{Schwartz's theorem for weak consistency}

Below, we state a result for posterior consistency
with respect to the weak topology due to \citet{schwartz1965} (see also \citet[Theorem 4.4.2]{ghosh2003}).
The result is a posterior consistency theorem for the density,
and thus relies on the assumption that the space of models $\FF$ is dominated by a $\sigma$-finite
measure $\mu$.

\begin{theorem}[Schwartz]
    \label{thm-schwartz}
    Let $\Pi$ be a prior on $\FF$
    and suppose $f_0$ is in the KL support of the prior $\Pi$.
    Then the posterior is \emph{weakly consistent} at $f_0$:
     i.e., for any weak neighborhood $U$ of $f_0$ the sequence of posterior distributions satisfies
\begin{align}
\Pi(U \,|\, X_{1:N}) \stackrel{N\rightarrow\infty}{\longrightarrow} 1, \quad
    \text{$f_0$-a.s.}
\end{align}
\end{theorem}

The above result assumes that the prior $\Pi$ is fixed.
Note that weak consistency also holds (with $f_0$-probability) for priors that
vary with $N$, provided that the sequence $\Pi_N$
satisfies the KL support condition
stated in \Cref{assumption-kl-sequence}
\citep[Theorem~6.25]{ghosal2017fundamentals}.

\section{Finite mixture models with an upper bound on the number of components}
\label{appendix-finite-prior}

In this section we consider a modification of the setting from the main paper
in which the prior $\Pi$ has support on only those finite mixtures with at
most $\tilde{k}$ components. We start by stating and proving our main result in this finite-support case.
Then we discuss why our conditions have changed slightly from \Cref{thm-inconsistency}.
Finally we demonstrate our finite-support theory in practice with an experiment.

\subsection{Result and proof}
Let $\FF(k)$ be the set of finite mixtures with exactly $k$ components for $k\leq \tilde{k}$.
We can apply the same proof technique in \Cref{sec-proof} to the present
case, provided that the mixture-density posterior concentrates on
weak neighborhoods of some compact subset of $\tilde{k}$-mixtures.

\begin{theorem}
\label{thm-finite-prior}
Suppose that the prior $\Pi$ has support on only those mixtures with at most $\tilde{k}$ components.
Assume that:
\begin{enumerate}
\vspace{-5pt}
\addtolength{\itemindent}{0cm}
    \item The posterior concentrates on weak neighborhoods of
        a weak-compact subset of $\FF(\tilde{k})$, and
    \item $\Psi$ is continuous, is mixture-identifiable, and has degenerate limits.
\vspace{-5pt}
\end{enumerate}
    Then the posterior on the number of components concentrates on ${\tilde{k}}$:
\begin{align}
    \Pi(\tilde{k} \,|\, X_{1:N}) \stackrel{N\rightarrow\infty}{\longrightarrow} 1
 \quad f_0\text{-a.s.}
\end{align}
\end{theorem}

\begin{proof}[Proof Sketch]
    By assumption, the posterior concentrates on weak neighborhoods of some weak-compact
    subset $\mathcal{A} \subseteq \FF(\tilde{k})$.
    It remains to show that there exists a weak neighborhood
    $U$ of $\mathcal{A}$ that, for all $k < \tilde k$, contains no $k$-mixtures of the family of $\Psi$.
    Suppose the contrary, i.e., that every such neighborhood contains a mixture of strictly less than $\tilde k$ components;
    then we can construct a sequence $(f_i)_{i=1}^\infty$ of mixtures of strictly less
    than $\tilde k$ components such that $f_i \Rightarrow \mathcal{A}$ (in the sense that the infimum of the weak metric
    between $f_i$ and elements of $\mathcal{A}$ converges to 0).
    Let $g_i$ be the corresponding sequence of mixing measures such that $f_i = F(g_i)$.
    Now we follow step 2 of the proof of the main theorem,
    with some slight modifications to account for the fact that $f_i$ converges weakly
    to a set rather than a single density.
    Suppose that $g_i(\Theta \setminus K) \to 0$ for some compact subset $K \subseteq \Theta$.
    Then following the proof of the main theorem, we have that $\FF_K$ and $\GG_K$ are compact,
    and so there is a weak-convergent subsequence of $F(\hat g_{i,K})$ that converges to some $f_0$; since $\mathcal{A}$
    is weak-closed, $f_0\in\mathcal{A}$. The remainder
    of this branch of the proof then follows the proof main theorem directly.
    Now for the other branch, suppose $g_i(\Theta \setminus K) \not\to 0$ for any compact $K\subseteq \Theta$.
    Then as in the main proof there is a sequence of parameters that is not relatively compact; so the
    corresponding sequence of components $\psi_i$ is either not tight or not $\mu$-wide.
    Since $\mathcal{A}$ is weak-compact, by Prokhorov's theorem $\mathcal{A}$ is tight, so
    $f_i$ must be tight, so $\psi_i$ must be tight. On the other hand, $\psi_i$ also must be $\mu$-wide,
    since otherwise replacing it with the singular sequence $\phi_i$ shows that $f_i$ would not converge
    weakly to $\mathcal{A}$. This concludes the second branch of the proof, and the result follows.
\end{proof}

\subsection{Discussion of the weak concentration condition}
Our main result in \Cref{thm-inconsistency} uses a KL support condition
to guarantee weak concentration of the posterior. In contrast, in \Cref{thm-finite-prior},
we do not impose any KL support condition and instead just directly
assume weak posterior concentration for the mixture density. First we discuss why this assumption remains
reasonable and then discuss why we chose to change the condition.

\paragraph{Reasonableness of the condition}
Note that the new weak-concentration assumption is actually weaker than the KL condition in the main paper---albeit potentially
substantially more difficult to verify. As a simple example of why this assumption is reasonable,
suppose we obtain data generated
from a Laplace distribution, and we use a mixture
model with Gaussian components and a prior that asserts that the mixture has at most 10 components.
Then we expect the posterior to concentrate on mixture densities that have exactly 10 components,
and in particular, the set of KL-closest mixture densities to the Laplace. Although many examples
will have a single closest such density, we state \Cref{thm-finite-prior} in such a way that
it allows for the case where the posterior concentrates on a \emph{compact set} of densities
(usually due to symmetry in the model).

\paragraph{Why change the condition}
In the main text, we assume---via the KL support condition, \Cref{assumption-kl-support}---that
the infimum of the KL divergence from the data generating distribution $f_0$ to
mixture distributions from the model is 0. In other words, we must be able to approximate $f_0$
arbitrarily well using mixture distributions from the model.
However, in the setting with a bounded number of components, this assumption typically does not hold. In
particular, the infimum KL
from the data-generating distribution $f_0$ to mixture distributions in the model is nonzero.
For example, in the previous Laplace versus Gaussian mixture example, we require an unbounded number
of components to achieve a vanishing KL divergence. If we are limited to 10 components, the infimum
KL is nonzero.

Demonstrating weak consistency with a reasonable amount of generality when the KL support
condition does not hold is challenging; see for instance,
\citet[Lemma~2.8]{kleijnphd} and \citet[Remark~4]{ramamoorthi2015posterior}.
Thus, we opt to require that weak concentration be verified directly
for each particular applied setting of interest, rather than attempting
to develop a general set of sufficient conditions.
The fact that we directly require weak convergence also means that we
do not need to make any stipulations about how data are generated. Therefore, in
contrast to the main theorem, we do not impose any such assumptions.

\begin{figure*}[t]
\centering
\begin{subfigure}[b]{0.46\linewidth}
    \centering
    \includegraphics[scale=0.44]{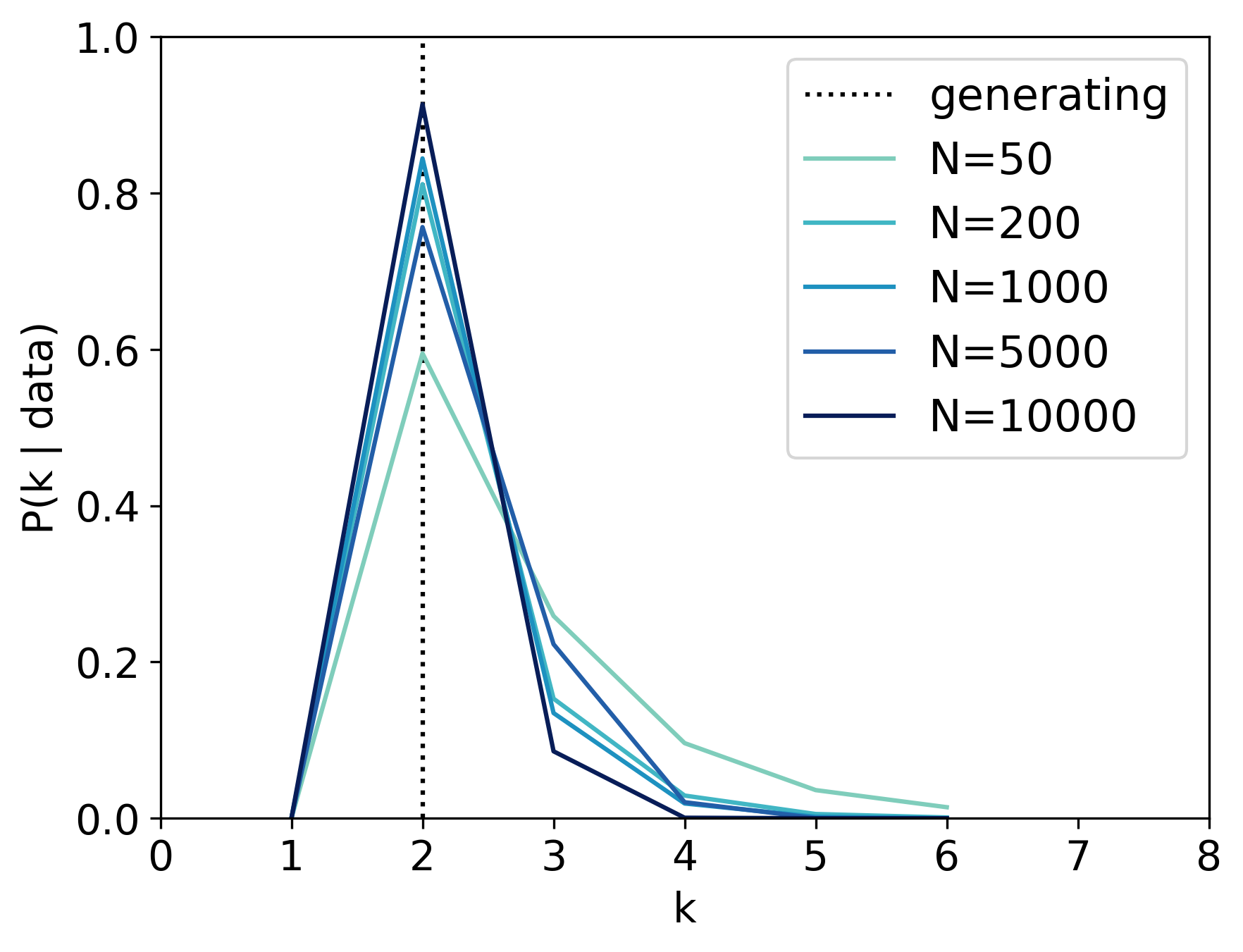}
    \caption{Gaussian mixture data}
    \label{fig-gaussian-finite}
\end{subfigure}
\begin{subfigure}[b]{0.46\linewidth}
    \centering
    \includegraphics[scale=0.44]{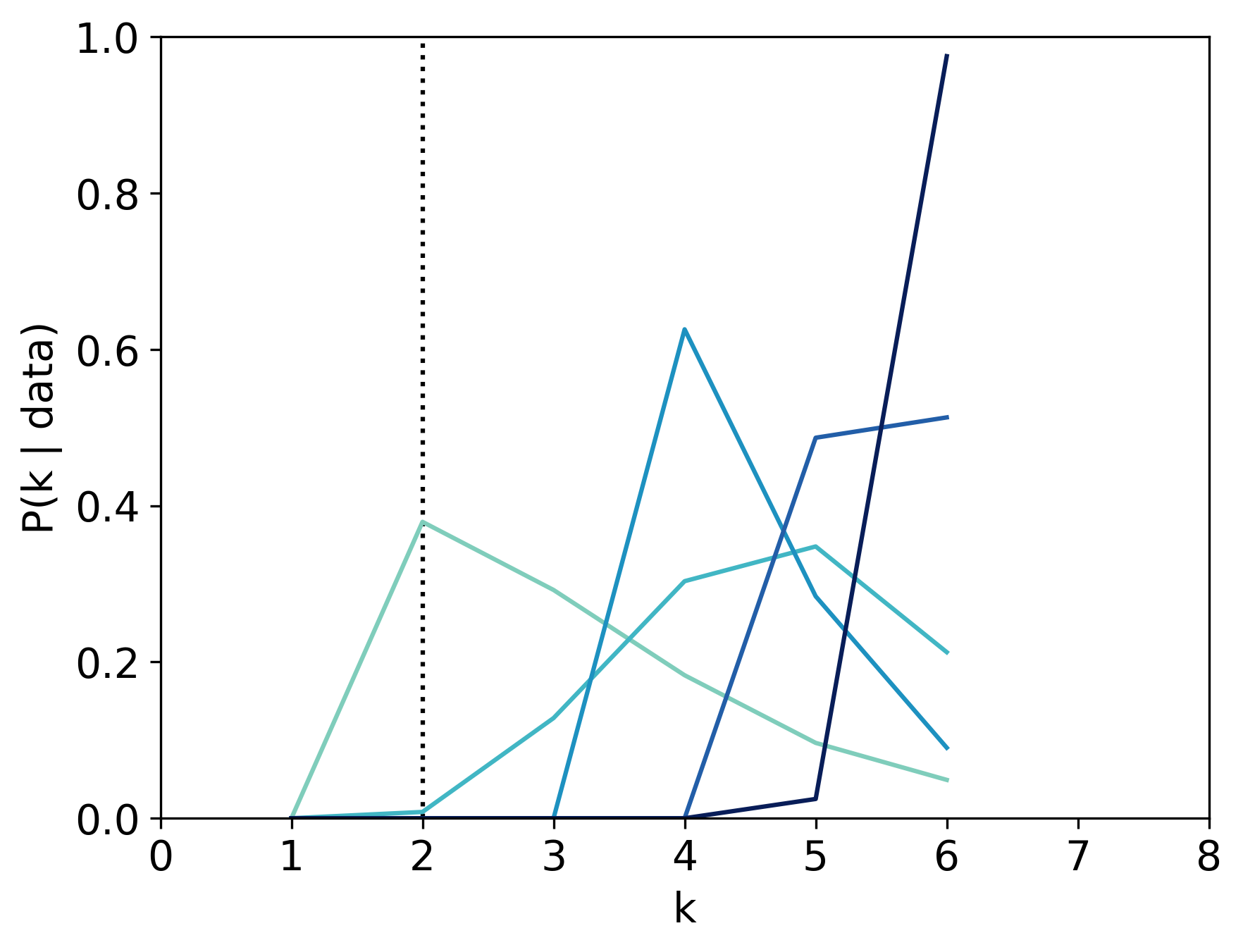}
    \caption{Laplace mixture data}
    \label{fig-laplace-finite}
\end{subfigure}
\caption{
    Well-specified and misspecified component families that use a prior with an
    upper bound on the number of components given by $k \sim \text{Unif}\{1,\ldots,6\}$.
    Posterior values for component counts $k$ with $k > 6$ are all zero, so we do not plot them.
}
    \label{figure-appendix-finite-prior}
\end{figure*}

\subsection{Experiments}

Now we demonstrate that the asymptotic behavior described by our theory occurs in practice.
In order to the study both the well-specified and misspecified cases,
we consider the same 2-component Gaussian and Laplace data described in
\Cref{ssec-synthetic}.
Here, we set the prior on the number of components to be a uniform distribution on $\{1,\ldots,6\}$.
The resulting posterior number of components appears in
 \Cref{figure-appendix-finite-prior}.
Here the well-specified model (Gaussian data) is consistent and concentrates
on the true generating number of components as $N$ grows
\citep{rousseau2011}.
On the other hand,
in the misspecified model (Laplace data), the posterior concentrates on
the largest possible number of components under the prior, in this case given by $\tilde{k}=6$.

\section{Proof of Proposition~2.2}
\label{appendix-gaussians}

Consider the multivariate Gaussian family
$\Psi = \left\{ \mathcal{N}(\nu, \Sigma) : \nu \in \reals^d, \, \Sigma\in \mathbb{S}^d_{++} \right\}$
with  parameter space
$\Theta = \reals^{d} \times \,\spd^d$, equipped with the topology induced by the Euclidean metric.
Let $(\lambda_j(\Sigma))_{j=1}^d$
denote the eigenvalues of the covariance matrix
$\Sigma \in \spd^d$ that satisfy
$\infty > \lambda_1(\Sigma) \geq \ldots \geq \lambda_d(\Sigma) > 0$.
Since the family of Gaussians is continuous and
mixture-identifiable \citep[Proposition~2]{yakowitz1968identifiability},
the main condition we need to verify is that the family has degenerate limits
(\Cref{defn-degenerate-limits}).
A useful fact is that if a sequence of Gaussian distributions is tight,
then the sequence of means and the eigenvalues of the covariance matrix is bounded.

\begin{lemma}
    \label{lemma-gaussians-tight}
    Let $(\psi_i)_{i \in \nats}$ be a sequence of Gaussian distributions with
    mean $\nu_i \in \reals^d$ and covariance $\Sigma_i \in \spd^d$.
    If $(\psi_i)_{i \in \nats}$ is a tight sequence of measures, then the sequences
    $(\nu_i)_{i \in \nats}$ and $(\lambda_1(\Sigma_i))_{i \in \nats}$ are bounded.
\end{lemma}

\begin{proof}
Let $Y_i$ denote a random variable with distribution $\psi_i$.
For each covariance matrix $\Sigma_i$, consider its eigenvalue decomposition
$\Sigma_i = U_i \Lambda_i U_i^\top$,
    where $U_i \in \reals^{d\times d}$ is an orthonormal matrix
    and $\Lambda_i \in \reals^{d\times d}$ is a diagonal matrix.
Then the random variable $Z_i = U_i^\top Y_i$ has distribution
$\mathcal{N}(U_i^\top \nu_i, \Lambda_i)$. If either $\|\nu_i\|_2 = \|U_i^\top \nu_i\|_2$
is unbounded or $\|\Lambda_i\|_F$ is unbounded, then $Z_i$ is not tight
 \citep[Example~25.10]{billingsley86}.
Since $Z_i$ and $Y_i$ lie in any ball centered at the origin with the same
probability, $Y_i$ is not tight.
\end{proof}

We now show that the multivariate Gaussian family
has degenerate limits.
\begin{proof}[Proof of \Cref{prop-gaussians}]
    If the parameters $(\theta_i)_{i \in \nats}$ are not a relatively compact
    subset of $\Theta$, then
    either some coordinate of the sequence
    of means $\nu_i$ diverges,
    $\lambda_1(\Sigma_i) \rightarrow \infty$,
    or
    $\lambda_d(\Sigma_i) \rightarrow 0$.
    If some coordinate of the mean $\nu_i$ diverges
    or the maximum eigenvalue diverges, i.e.,
    $\lambda_1(\Sigma_i) \rightarrow \infty$,
    then the sequence $(\psi_{\theta_i})$ is not tight
    by \Cref{lemma-gaussians-tight}.
    On the other hand, if $\lambda_d(\Sigma_i) \rightarrow 0$ as $i \rightarrow \infty$,
    then $\psi_{\theta_i}$ converges weakly to a sequence of degenerate Gaussian measures
    that concentrate on $C_i = \left\{x \in \reals^d : (x-\nu_i)^\top u_{d,i} = 0\right\}$,
    where $u_{d,i}$ is the $d^\text{th}$ eigenvector of $\Sigma_i$.
    Note that $\mu(C_i) = 0$ for Lebesgue measure $\mu$;
    so if we define $C = \cup_i C_i$ in the setting of \Cref{defn-muwide}, the sequence
    is not $\mu$-wide.
\end{proof}

We can generalize \Cref{prop-gaussians} beyond multivariate Gaussians
to mixture-identifiable location-scale families, as shown in \Cref{prop-locscale}.
Examples of such families include
the multivariate Gaussian family,
the Cauchy family, the logistic family, the von Mises family, and
generalized extreme value families.
The proof is similar to that of \Cref{prop-gaussians}.
\begin{proposition}\label{prop-locscale}
Suppose $\Psi$ is a location-scale family
that is mixture-identifiable and absolutely continuous with
respect to Lebesgue measure $\mu$, i.e.,
\[
    \frac{d\Psi}{d\mu} =
    \left\{
        |\Sigma|^{-\nicefrac{1}{2}} \, \varphi\left(\Sigma^{-\nicefrac{1}{2}}(x - \nu)\right):
    \nu \in \reals^d, \Sigma \in \spd^d
    \right\},
\]
where $\varphi: \reals^d \rightarrow \reals$ is a probability density function.
Then $\Psi$ satisfies \Cref{assumption-regularity}.
\end{proposition}

\section{Additional related work}

Priors for microclustering behavior  have been a recent focus in the
Bayesian nonparametrics literature \citep{zanella2016flexible,klami2016probabilistic}.
Since having a fixed number of components across dataset sizes $N$ would be incompatible with
sublinear growth (in $N$) of cluster size across all clusters,
we expect divergence issues similar in flavor to those in
\citet{miller13,miller2014}.

\bibliographystyle{plainnat-mod}
\bibliography{main}

\end{document}